\PassOptionsToPackage{unicode}{hyperref}
\PassOptionsToPackage{hyphens}{url}
\PassOptionsToPackage{dvipsnames,svgnames,x11names}{xcolor}
\documentclass[
  12pt,
  a4paper,
]{article}
\usepackage{xcolor}
\usepackage[margin=1in]{geometry}
\usepackage{amsmath,amssymb}
\usepackage{iftex}
\ifPDFTeX
  \usepackage[T1]{fontenc}
  \usepackage[utf8]{inputenc}
  \usepackage{textcomp} 
\else 
  \usepackage{unicode-math} 
  \defaultfontfeatures{Scale=MatchLowercase}
  \defaultfontfeatures[\rmfamily]{Ligatures=TeX,Scale=1}
\fi
\usepackage{lmodern}
\ifPDFTeX\else
\fi
\IfFileExists{upquote.sty}{\usepackage{upquote}}{}
\IfFileExists{microtype.sty}{
  \usepackage[]{microtype}
  \UseMicrotypeSet[protrusion]{basicmath} 
}{}
\usepackage{setspace}
\makeatletter
\@ifundefined{KOMAClassName}{
  \IfFileExists{parskip.sty}{%
    \usepackage{parskip}
  }{
    \setlength{\parindent}{0pt}
    \setlength{\parskip}{6pt plus 2pt minus 1pt}}
}{
  \KOMAoptions{parskip=half}}
\makeatother
\makeatletter
\ifx\paragraph\undefined\else
  \let\oldparagraph\paragraph
  \renewcommand{\paragraph}{
    \@ifstar
      \xxxParagraphStar
      \xxxParagraphNoStar
  }
  \newcommand{\xxxParagraphStar}[1]{\oldparagraph*{#1}\mbox{}}
  \newcommand{\xxxParagraphNoStar}[1]{\oldparagraph{#1}\mbox{}}
\fi
\ifx\subparagraph\undefined\else
  \let\oldsubparagraph\subparagraph
  \renewcommand{\subparagraph}{
    \@ifstar
      \xxxSubParagraphStar
      \xxxSubParagraphNoStar
  }
  \newcommand{\xxxSubParagraphStar}[1]{\oldsubparagraph*{#1}\mbox{}}
  \newcommand{\xxxSubParagraphNoStar}[1]{\oldsubparagraph{#1}\mbox{}}
\fi
\makeatother

\usepackage{longtable,booktabs,array}
\usepackage{calc} 
\usepackage{etoolbox}
\makeatletter
\patchcmd\longtable{\par}{\if@noskipsec\mbox{}\fi\par}{}{}
\makeatother
\IfFileExists{footnotehyper.sty}{\usepackage{footnotehyper}}{\usepackage{footnote}}
\makesavenoteenv{longtable}
\usepackage{graphicx}
\makeatletter
\newsavebox\pandoc@box
\newcommand*\pandocbounded[1]{
  \sbox\pandoc@box{#1}%
  \Gscale@div\@tempa{\textheight}{\dimexpr\ht\pandoc@box+\dp\pandoc@box\relax}%
  \Gscale@div\@tempb{\linewidth}{\wd\pandoc@box}%
  \ifdim\@tempb\p@<\@tempa\p@\let\@tempa\@tempb\fi
  \ifdim\@tempa\p@<\p@\scalebox{\@tempa}{\usebox\pandoc@box}%
  \else\usebox{\pandoc@box}%
  \fi%
}
\def\fps@figure{htbp}
\makeatother

\usepackage[]{natbib}
\usepackage[utf8]{inputenc}
\usepackage{mathtools}
\usepackage{amsthm}
\usepackage{algorithm}
\usepackage{algpseudocode} 
\usepackage{tikz}
\usepackage{comment}
\usepackage{authblk}
\usepackage{bm}

\newcommand{\EE}{\mathbb{E}}
\newcommand{\PP}{\mathbb{P}}

\newcommand{\fdr}{\operatorname{FDR}}
\newcommand{\fdp}{\operatorname{FDP}}

\newcommand{\cA}{\mathcal{A}}

\newcommand{\cC}{\mathcal{C}}

\newcommand{\cG}{\mathcal{G}}

\newcommand{\cK}{\mathcal{K}}

\newcommand{\cN}{\mathcal{N}}

\newcommand{\cP}{\mathcal{P}}

\newcommand{\cV}{\mathcal{V}}

\newcommand{\ebh}{\textnormal{e-BH}}

\newcommand{\Uebh}{\textnormal{U-eBH}}
\newcommand{\Rone}{\textnormal{R1-eBH}}
\newcommand{\Rtwo}{\textnormal{R2-eBH}}
\newcommand{\Rboth}{\textnormal{R-eBH}}
\newcommand{\Jebh}{\textnormal{J-eBH}}
\newcommand{\BY}{\textnormal{BY}}
\newcommand{\BH}{\textnormal{BH}}
\newcommand{\Uby}{\textnormal{U-BY}}
\newcommand{\DstarBY}{\ensuremath{D^{\star}_{\textnormal{BY}}}}

\newtheorem{lemma}{Lemma}
\makeatletter
\@ifpackageloaded{caption}{}{\usepackage{caption}}
\AtBeginDocument{%
\ifdefined\contentsname
  \renewcommand*\contentsname{Table of contents}
\else
  \newcommand\contentsname{Table of contents}
\fi
\ifdefined\listfigurename
  \renewcommand*\listfigurename{List of Figures}
\else
  \newcommand\listfigurename{List of Figures}
\fi
\ifdefined\listtablename
  \renewcommand*\listtablename{List of Tables}
\else
  \newcommand\listtablename{List of Tables}
\fi
\ifdefined\figurename
  \renewcommand*\figurename{Figure}
\else
  \newcommand\figurename{Figure}
\fi
\ifdefined\tablename
  \renewcommand*\tablename{Table}
\else
  \newcommand\tablename{Table}
\fi
}
\@ifpackageloaded{float}{}{\usepackage{float}}
\floatstyle{ruled}
\@ifundefined{c@chapter}{\newfloat{codelisting}{h}{lop}}{\newfloat{codelisting}{h}{lop}[chapter]}
\floatname{codelisting}{Listing}

\usepackage{amsthm}
\theoremstyle{plain}
\newtheorem{theorem}{Theorem}[section]
\theoremstyle{plain}
\newtheorem{proposition}{Proposition}[section]
\theoremstyle{definition}
\newtheorem{definition}{Definition}[section]
\theoremstyle{remark}
\AtBeginDocument{}

\newtheorem{corollary}{Corollary}[section]

\makeatother
\makeatletter
\@ifpackageloaded{caption}{}{\usepackage{caption}}
\@ifpackageloaded{subcaption}{}{\usepackage{subcaption}}
\makeatother
\usepackage{bookmark}
\IfFileExists{xurl.sty}{\usepackage{xurl}}{} 
\hypersetup{
  pdftitle={Title},
  pdfauthor={Bose, Bartroff},
  colorlinks=true,
  linkcolor={blue},
  filecolor={Maroon},
  citecolor={Blue},
  urlcolor={Blue},
  pdfcreator={LaTeX}}

\allowdisplaybreaks

\title{Optimal Certified Multiple Testing with E-values}

\author{Soumyabrata Bose\thanks{\texttt{sbose@utexas.edu}} }
\author{Jay Bartroff\thanks{\texttt{bartroff@austin.utexas.edu}}}
\affil{
    Department of Statistics and Data Sciences \\
    University of Texas at Austin
}

\date{\today}
\begin{document}
\maketitle
\begin{abstract}
We study randomized multiple testing with e-values under arbitrary dependence and introduce a conditional certification constraint: Given the observed evidence, each coordinate's expected rejection certificate must not exceed its e-value. Existing randomized e-BH improvements satisfy this constraint but have no established optimality guarantee. We prove that the maximum conditional inclusion probability for hypothesis $i$, given e-values~$x$, is $1 \wedge (\alpha / K) \sum_j \min(x_i, x_j)$, and construct one rule attaining these bounds simultaneously. The rule uses a single independent uniform variable, rejects a prefix of the ordered e-values, and has a unique conditional rejection-set law among rules attaining all the bounds. It therefore maximizes every nonnegative weighted sum of inclusion probabilities, including expected true discoveries under any truth configuration. Under a shared uniform draw, its rejection set contains those of e-BH and U-eBH. The construction extends to BY-calibrated p-values, and simulations under positive and negative dependence show power gains over existing randomized procedures.
\end{abstract}

\setstretch{2}

\section{Introduction} \label{sec-intro}

Multiple testing problems routinely involve dependence among the hypotheses that is unknown or too complex to model, so a procedure must control the false discovery rate (FDR) without relying on a particular dependence structure. The Benjamini--Hochberg (BH) procedure controls the FDR when the null p-values are independent \citep{Benjamini95} and remains valid under positive regression dependence on a subset \citep{Benjamini01}; under arbitrary dependence, control is recovered by a correction that lowers the level at which each hypothesis is tested, at a cost in power that grows with the number of hypotheses \citep{Benjamini01}. The e-BH procedure of \citet{Wang22} avoids the correction by operating on e-values, nonnegative statistics whose expectations under the null hypotheses are bounded by one \citep{Vovk21}. Each e-value is only required to satisfy its own marginal expectation bound, and FDR control follows under arbitrary dependence among the coordinates. Dependence-robust validity therefore follows from the marginal bounds alone, and the question of how effectively the observed evidence is converted into discoveries is left open by the validity guarantee itself.

Improving the power of e-value-based multiple testing while preserving its dependence-robust guarantees has consequently attracted considerable attention. One established route exploits information beyond the marginal evidence. Weighting assigns hypothesis-specific weights to p-values before applying a multiple testing procedure, and increases discoveries when the weights favor the alternatives, with validity preserved as long as the weights are independent of the p-values \citep{Genovese06}; \citet{Ignatiadis24} extend the construction to e-values, using an e-value as an unnormalized weight for a p-value. A second route uses distributional information about the null. \citet{Wang22} show that the power of e-BH can be boosted when the marginal null distribution of each e-value is known, \citet{Lee24} construct boosted e-values through conditional calibration, and \citet{Blier26} derive improved e-value thresholds under additional distributional assumptions. Both routes depend on the accuracy of the auxiliary or distributional information, and their validity is tied to the assumptions under which that information is obtained.

A different route, introduced by \citet{Xu25}, uses external randomization and requires no additional information. The underlying operation is stochastic rounding, which replaces a nonnegative evidence value by one of the two adjacent grid points with probabilities chosen so that its expectation is preserved. The operation predates its application to multiple testing: it is standard in numerical analysis, where it controls accumulated rounding error \citep{Croci22}, and for a single e-value at a fixed threshold it is the randomized Markov inequality of \citet{Ramdas26b}, which rejects with probability equal to the product of the test level and the e-value, truncated at one, and is never less powerful than the deterministic threshold rule. Extending the operation to the data-dependent thresholds of e-BH, and to jointly randomized decisions across coordinates, is due to \citet{Xu25}. Their procedures coordinate the rounding across hypotheses through either a shared uniform variable or coordinate-specific uniform variables, and each preserves FDR control under arbitrary dependence while producing at least as many rejections as e-BH. The coupling matters because e-BH uses thresholds determined by the entire evidence vector. Randomizations preserving the same coordinatewise expectations can therefore produce different rejection sets and different power.

The joint law of the randomization is thus a free parameter within the construction, and the choice among its realizations is left open. At a fixed evidence vector, a rule of this form retains the e-BH rejection set and adds a random set supported on the remaining coordinates, subject to a conditional expectation constraint for each coordinate not already rejected. \citet{Xu25} state the resulting optimization problem but do not solve it; in their conclusion they observe that no optimality has been established for any of their procedures, and they leave both the identification of an optimal scheme and the existence of an efficient algorithm open. Optimization of a discovery objective under error constraints has an established literature of its own, in which \citet{Rosset22} characterize optimal and maximin rejection policies for specified testing models, but that tradition optimizes over the sampling distribution of the data. Here the statistical model is fixed and the data are given; the free variable is the joint law of the external randomization, and the quantity to be maximized is conditional on the observed evidence. The formulation of \citet{Xu25} does not characterize the inclusion probabilities attainable under the conditional expectation constraints, and it does not determine whether a single randomized rule can attain the pointwise maximum of those probabilities at every hypothesis simultaneously. These are distinct requirements, since a rule can maximize an aggregate discovery objective without attaining the pointwise maximum at every hypothesis. Three questions follow: 
\begin{itemize}
    \item What limits do the conditional expectation constraints impose on the inclusion probability of each hypothesis? 
    \item Can a single randomized rule attain all of these limits simultaneously? 
    \item If such a rule exists, is its conditional rejection-set law uniquely determined?
\end{itemize}

This paper addresses these questions through a framework of \emph{conditional certification}. We associate with each realized rejection set a vector of certificates determined by its cardinality, motivated by the self-consistency property underlying e-BH \citep{Blanchard08}, and require that the conditional expectation of each certificate, given the observed evidence, be bounded by the corresponding e-value. Because the requirement constrains each coordinate separately at a fixed evidence vector, any bound it implies applies to each hypothesis on its own and holds under every distribution compatible with the e-value condition. Joint rounding is one way to satisfy these constraints, but it is not the only one, and it is not the one adopted here: certification imposes them on the distribution of the rejection set directly, so the rules it admits need not be obtainable by rounding the evidence, and the class it defines is correspondingly larger than the class of joint roundings. The construction of \citet{Xu25} motivates the optimization problem. The certification framework formulates that problem over the larger class of rules it admits, and the analysis in this paper solves it explicitly.

We derive sharp upper bounds on the conditional probability that each hypothesis can be rejected under certification, and construct a single randomized procedure that attains all of them simultaneously. The rule uses one independent uniform random variable and rejects a prefix of the ordered e-values, a form that follows from the analysis rather than being imposed at the outset. Any rule maximizing a weighted sum of inclusion probabilities with strictly positive weights has the same conditional rejection-set law.
The procedure therefore maximizes every nonnegative weighted sum of conditional inclusion probabilities, including the expected number of discoveries and the expected number of true discoveries under any configuration of true and false nulls, and the optimality is uniform over the certified class rather than relative to a particular competitor. We also establish a pathwise nesting: for every evidence vector and every realization of the shared uniform variable, the e-BH rejection set is contained in that of the \textnormal{U-eBH} procedure of \citet{Xu25}, which is contained in that of the proposed rule.

The optimality is relative to the class of conditionally certified rules, not to all FDR-controlling randomized procedures. Compound e-values \citep{Ignatiadis24} and e-value-based closure principles \citep{Xu25b} fall outside this class, and the results do not constrain them. Within the class, the characterization is complete: the bounds are sharp, and a rule attaining them is exhibited. Numerical experiments under positive and negative dependence show that the theoretical gains translate into power improvements over the existing randomized e-BH procedures, largest at intermediate signal strengths, with attained FDR below the target throughout.

The remainder of the paper is organized as follows. Section~\ref{sec-setup} introduces the testing framework and notation. Section~\ref{sec-certification} develops conditional certification and establishes the bounds on rejection probabilities, and Section~\ref{sec-optimality} constructs the optimal procedure and proves its optimality and uniqueness properties. Section~\ref{sec-simulation} presents numerical comparisons with the existing e-BH procedures, and the appendix contains the technical verifications and the descriptions of the competing randomized procedures.

\section{Setup and Notation}\label{sec-setup}

We consider the problem of simultaneously testing $K$ hypotheses using e-values under arbitrary dependence, with target false discovery rate $\alpha \in (0, 1)$. We allow the decision rules to depend on an external source of randomization independent of the evidence. We first define the evidence model and the randomized decision rules, then the performance criteria and the e-BH benchmark.

\subsection{The Testing Problem and Randomized Rules} \label{sec-setup-evidence}

Let $\cP$ be a class of distributions on a common measurable space. For each $i \in [K]$, let $H_i \subseteq \cP$ denote the $i$\textsuperscript{th} null hypothesis. For $P \in \cP$, let $\cN(P) = \{i \in [K] : P \in H_i\}$ denote the set of true null hypotheses, with $K_0(P) = |\cN(P)|$ true nulls and $K_1(P) = K - K_0(P)$ alternatives. We write $\cN^c = [K] \setminus \cN$ for the complement. We observe an evidence vector $X = (X_1, X_2, \dots, X_K) \in [0, \infty)^K$, finite almost surely. Each coordinate is an e-value for its hypothesis, i.e., 
\begin{equation*}
\EE_P[X_i] \le 1
\:\: \text{for every } P \in H_i.
\end{equation*}
The coordinates may be arbitrarily dependent, and no expectation restriction is imposed on $X_i$ when $P \notin H_i$.

Let $V$ be an auxiliary random element with distribution $Q$ on a measurable space $\cV$, independent of $X$ under every $P \in \cP$. A randomized rejection rule is a measurable map $D : [0, \infty)^K \times \cV \to 2^{[K]}$. At evidence $X$, the reported rejection set is $D(X, V)$. A deterministic rule is a special case in which $D(x, V)$ does not depend on $V$. We use $\PP_V$ and $\EE_V$ for probability and expectation over the randomization at fixed $x$, and $\PP_{P, V}$ and $\EE_{P, V}$ when both the data and the randomization are integrated. The conditional inclusion probability of hypothesis $i$ under rule $D$ is
$$
\pi_i^D(x) = \PP_V \{i \in D(x, V)\}.
$$
The rejection count is $R_D(x, V) = |D(x, V)|$. 

\subsection{Error Control and Discovery Objectives} \label{sec-setup-error}

For a realized rejection set $A \subseteq [K]$ and a distribution $P \in \cP$, the false discovery proportion is
$$
\fdp_P(A) = \frac{|A \cap \cN(P)|}{|A| \vee 1}.
$$
For a randomized rule $D$, the false discovery rate is $\fdr_P(D) = \EE_{P, V} \left[\fdp_P(D(X, V)) \right]$. The rule $D$ is valid at level $\alpha$ if $\fdr_P(D) \le \alpha$ for every $P \in \cP$. When $K_1(P) > 0$, the average power of $D$ under $P$ is
$$
\operatorname{Power}_P(D) = \frac{\EE_{P, V} \left[|D(X, V) \cap \cN(P)^c| \right]}{K_1(P)} = 
\frac{1}{K_1(P)} \sum_{i \notin \mathcal N(P)} \mathbb E_P [\pi_i^D(X)]. 
$$

We compare conditional inclusion probabilities $\pi_i^D (x)$ at fixed inputs in the following sections. Averaging over $X$ then yields comparisons of average power.

\subsection{The e-BH Rule and Self-Consistency}\label{sec-setup-ebh}

For a fixed $x$, we write $x_{[1]} \ge x_{[2]} \ge \cdots \ge x_{[K]}$ for its coordinates in decreasing order. Define the e-BH thresholds and selected rank by
$$
q_r = \frac{K}{\alpha r},
\:\:
r_{\operatorname e} (x) = \max \left(\{r \in [K] : x_{[r]} \ge q_r\} \cup \{0\} \right).
$$
The e-BH rejection set is
\begin{equation} \label{def-ebh}
D_{\operatorname e} (x) =
\begin{cases}
\{i \in [K] : x_i \ge q_{r_{\operatorname e} (x)}\} & r_{\operatorname e} (x) > 0, \\
\emptyset & r_{\operatorname e} (x) = 0.
\end{cases}
\end{equation}
By construction, $|D_{\operatorname e} (x)| = r_{\operatorname e} (x)$. Under the e-value condition, the e-BH procedure controls the false discovery rate \citep{Wang22} with
$$
\fdr_P(D_{\operatorname e}) \le \frac{\alpha K_0(P)}{K} \le \alpha
\:\: \text{for every } P \in \cP.
$$

\begin{definition} \label{def-self-consistency}
A rejection set $A \subseteq [K]$ is \emph{self-consistent} with respect to $x$ at level $\alpha$ if
$$
i \in A
\:\: \Longrightarrow \:\:
x_i \ge \frac{K}{\alpha(|A| \vee 1)}.
$$
The empty set is self-consistent with respect to every $x$.
\end{definition}

For every $x$, the e-BH rejection set~$D_{\operatorname e}(x)$ is self-consistent with respect to $x$, and every rejection set that is self-consistent with respect to $x$ is contained in $D_{\operatorname e}(x)$. The e-BH rule therefore returns the maximal self-consistent set at every input, and it is the benchmark against which the procedures developed below are measured.

\section{Conditional Certification and Inclusion Bounds} \label{sec-certification}

The self-consistency condition of Section~\ref{sec-setup-ebh} requires each rejected coordinate to clear an evidence threshold determined by the rejection count. We express this requirement through \textit{certificates} and study the resulting inequalities in expectation over the external randomization. The development proceeds in three stages: we first impose the inequalities on true-null coordinates and establish FDR control, then derive a bound on conditional inclusion probabilities under that constraint, and finally introduce a condition that can be verified without identifying the true-null set.

\subsection{Certificates and FDR Control} \label{sec-certification-fdr}

For a realized rejection set $A$, define its \emph{certificate} $Z(A) = (Z_1(A), Z_2(A), \dots, Z_K(A))$ by
$$
Z_i(A) = \frac{K}{\alpha (|A| \vee 1)} \mathbf 1 \{i \in A\}, \:\: i \in [K].
$$
The certificate $Z_i(A)$ equals $K / (\alpha |A|)$ when $i \in A$ and zero otherwise. The self-consistency condition of Definition~\ref{def-self-consistency} is equivalent to the pointwise inequality $Z_i(A) \le x_i$ for every  $i \in [K]$. For $i \notin A$, $Z_i(A) = 0 \le x_i$ holds automatically. 

For every $P \in \cP$, the false discovery proportion can be written in terms of the certificates of the true nulls. Substituting the definition of $Z_i(A)$ gives
$$
\fdp_P(A) = \frac{|A \cap \cN(P)|}{|A| \vee 1} = \frac{\alpha}{K} \sum_{i \in \cN (P)} \frac{K}{\alpha (|A| \vee 1)} \mathbf 1 \{i \in A\} = \frac{\alpha}{K} \sum_{i \in \cN (P)} Z_i(A). 
$$
Since only the certificates of the true nulls appear, we impose the conditional requirement on those coordinates. 

Fix a set $\cN \subseteq [K]$. A rule $D$ satisfies \emph{null-only conditional certification} at $\cN$ if
\begin{equation}\label{eq-null-only}
\EE_V[Z_i(D(x, V))] \le x_i, \:\: \text{for every } x \in [0, \infty)^K \text{ and } i \in \cN,
\end{equation}
where the expectation is over $V$ at fixed $x$. We refer to the right-hand side $x_i$ as the certificate budget for coordinate $i$. The certification condition \eqref{eq-null-only} is stated in terms of the certificate $Z(D)$. We now characterize it by the existence of a nonnegative evidence vector that supports the rejection set through self-consistency and satisfies the same conditional budgets.

A measurable vector $S : [0, \infty)^K \times \cV \to [0, \infty)^K$ \emph{supports} a rule $D$ if $D(x, V)$ is self-consistent with respect to $S(x, V)$ for every fixed $x$, i.e.,
$$
i \in D(x, V) \:\: \Longrightarrow \:\: S_i(x, V) \ge \frac{K}{\alpha (|D(x, V)| \vee 1)}. 
$$

\begin{proposition} \label{prop-characterization}
Fix a set $\cN \subseteq [K]$. A rule $D$ satisfies \eqref{eq-null-only} at $\cN$ if and only if there exists a nonnegative supporting vector $S$ such that
$$
\EE_V[S_i(x, V)] \le x_i, \:\: \text{for every } x \in [0, \infty)^K \text{ and } i \in \cN.
$$
\end{proposition}

\begin{proof}
Suppose first that $D$ admits a supporting vector $S$ with the stated budgets. For every $(x, V)$ and every $i$, self-consistency and nonnegativity give
$$
Z_i(D(x, V)) \le S_i(x, V), 
$$
as $Z_i = 0 \le S_i$ whenever $i \notin D(x, V)$, and $Z_i = K/(\alpha |D(x, V)|) \le S_i$ whenever $i \in D(x, V)$. Taking expectations and applying the budget on $S_i$ yields
$$
\EE_V[Z_i(D(x, V))] \le \EE_V[S_i(x, V)] \le x_i, \:\: i \in \cN.
$$

Conversely, if \eqref{eq-null-only} holds, take $S_i(x, V) = Z_i(D(x, V))$. The vector $S$ supports $D$ by construction, since $S_i(x, V) = K/(\alpha |D(x, V)|)$ on the rejection set, and the inequalities $\EE_V[S_i(x, V)] \le x_i$ required by the proposition are exactly \eqref{eq-null-only}.
\end{proof}

The certificate is the coordinatewise smallest nonnegative vector supporting a given rejection set. Imposing conditional budgets on $Z(D)$ is therefore equivalent to requiring the existence of supporting evidence satisfying those budgets.

\begin{proposition} \label{prop-fdr-null}
Fix $P \in \cP$. If a rule $D$ satisfies \eqref{eq-null-only} at $\cN(P)$, then $\fdr_P(D) \le {\alpha K_0(P)} / {K} \le \alpha$. 
\end{proposition}

\begin{proof}
By the certificate representation of the false discovery proportion and iterated expectation, we have 
$$
\begin{aligned}
\fdr_P(D) 
=&\: \frac{\alpha}{K} \sum_{i \in \cN(P)} \EE_{P, V}[Z_i(D(X, V))] \\
=&\: \frac{\alpha}{K} \sum_{i \in \cN(P)} \EE_P\!\left[\EE_V[Z_i(D(X, V))]\right] \\
\le&\: \frac{\alpha}{K} \sum_{i \in \cN(P)} \EE_P[X_i] \\
\le&\: \frac{\alpha K_0(P)}{K}.
\end{aligned}
$$
\end{proof}
Null-only conditional certification is sufficient for FDR control. We now derive the inclusion bound that this condition implies.

\subsection{An Oracle Inclusion Bound} \label{sec-certification-oracle}

Fix a deterministic input $x \in [0, \infty)^K$, a set $\cN \subseteq [K]$, and a rule $D$ satisfying \eqref{eq-null-only} at $\cN$. We ask what these inequalities imply for the conditional inclusion probability $\pi_i^D(x) = \PP_V\{i \in D(x, V)\}$. The calculation is pointwise in $x$ and requires no distributional assumptions on the evidence. The starting point is a certificate identity that expresses the selection of a coordinate through its certificate and the certificates of the remaining coordinates. For every realized set $A$ and every coordinate $i$, we have
\begin{equation}\label{eq-cert-identity}
\bm 1\{i \in A\} = \frac{\alpha}{K} \sum_{j = 1}^K \min\{Z_i(A), Z_j(A)\}. 
\end{equation}
To see this, if $i \notin A$, then $Z_i(A) = 0$, so every term in the sum vanishes and both sides equal $0$. If $i \in A$, then $Z_i(A) = K/(\alpha |A|)$. The terms with $j \in A$ equal $K / (\alpha |A|)$, and the terms with $j \notin A$ vanish. Since there are $|A|$ terms of the first kind, the sum equals $K / \alpha$, and multiplying by $\alpha / K$ gives $1$. We now apply \eqref{eq-cert-identity} to bound the conditional inclusion probability of a constrained coordinate.

\begin{proposition} \label{prop-oracle-bound}
Let $D$ satisfy \eqref{eq-null-only} at $\cN$, and let $K_1 = K - |\cN|$. Then for every $x \in [0, \infty)^K$ and every $i \in \cN$,
$$
\pi_i^D(x) \le 1 \wedge \frac{\alpha}{K} \left[ K_1 x_i + \sum_{j \in \cN} \min(x_i, x_j) \right].
$$
\end{proposition}

\begin{proof}
Fix $i \in \cN$. Taking expectations in \eqref{eq-cert-identity} with $A = D(x, V)$ and writing $Z_j = Z_j(D(x, V))$ gives
$$
\pi_i^D(x) = \frac{\alpha}{K} \left[\sum_{j \in \cN} \EE_V[\min(Z_i, Z_j)] + \sum_{j \notin \cN} \EE_V[\min(Z_i, Z_j)] \right].
$$
For $j \in \cN$, both certificates satisfy \eqref{eq-null-only}, and since the minimum is bounded above by either argument,
$$
\EE_V[\min(Z_i, Z_j)] \le \min\{\EE_V[Z_i], \EE_V[Z_j]\} \le \min(x_i, x_j).
$$
For $j \notin \cN$, only the certificate of $i$ is constrained, so
$$
\EE_V[\min(Z_i, Z_j)] \le \EE_V[Z_i] \le x_i.
$$
There are $K_1 = K - |\cN|$ coordinates with $j \notin \cN$. Substituting these bounds gives
$$
\pi_i^D(x) \le \frac{\alpha}{K} \left[\sum_{j \in \cN} \min(x_i, x_j) + K_1 x_i \right].
$$
Combining with the elementary bound $\pi_i^D(x) \le 1$ gives the stated inequality.
\end{proof}

This is an oracle bound: it depends on the true configuration through $\cN$ and $K_1$. The argument requires $i \in \cN$, since it uses the certificate budget for coordinate $i$; for $i \notin \cN$, null-only certification gives no improvement over the elementary bound $\pi_i^D(x) \le 1$. Since $\cN$ is unknown, we next impose the certification condition on every coordinate.

\subsection{Certification on All Coordinates} \label{sec-certification-class}

We impose the certificate condition on every coordinate and refer to the resulting rules as conditionally certified. The following sections study optimality within this class.
\begin{definition} \label{def-class-C}
The \emph{conditionally certified} class $\cC$ consists of all randomized rules $D$ satisfying
\begin{equation}\label{eq-all-coord}
\EE_V[Z_i(D(x, V))] \le x_i, \:\: \text{for every } x \in [0, \infty)^K \text{ and every } i \in [K].
\end{equation}
\end{definition}

Every $D \in \cC$ satisfies $\fdr_P(D) \le {\alpha K_0(P)} / {K}$ for every $P \in \cP$, by Proposition~\ref{prop-fdr-null} applied with $\cN = \cN(P)$. By Proposition~\ref{prop-characterization} applied with $\cN = [K]$, a rule $D$ belongs to $\cC$ if and only if it admits a nonnegative supporting vector $S$ with
$$
\EE_V[S_i(x, V)] \le x_i, \:\: \text{for every } x \in [0, \infty)^K \text{ and every } i \in [K].
$$
If $S$ satisfies the inequalities and $D(x, V) = D_{\operatorname e}(S(x, V))$, then $D \in \cC$, since e-BH is self-consistent with respect to its input. For an arbitrary supporting vector, self-consistency only implies $D(x, V) \subseteq D_{\operatorname e}(S(x, V))$.

The ordinary e-BH rule belongs to $\cC$: its self-consistency gives the pointwise inequalities $Z_i(D_{\operatorname e}(x)) \le x_i$ for every $i$, which imply \eqref{eq-all-coord} after taking expectations. Detailed membership arguments for the procedures of \citet{Xu25} are deferred to the Appendix.

Applying Proposition~\ref{prop-oracle-bound} with $\cN = [K]$ gives, for every $D \in \cC$, every $x \in [0, \infty)^K$, and every $i \in [K]$, $\pi_i^D(x) \le 1 \wedge G_i(x)$, where 
$$
G_i(x) = \frac{\alpha}{K} \sum_{j = 1}^K \min(x_i, x_j).
$$

\section{Optimality under Conditional Certification} \label{sec-optimality}

Section~\ref{sec-certification-class} established that every $D \in \cC$ satisfies the coordinatewise bound
$$
\pi_i^D(x) \le 1 \wedge G_i(x), \:\: G_i(x) = \frac{\alpha}{K} \sum_{j = 1}^K \min(x_i, x_j),
$$
for every $x \in [0, \infty)^K$ and every $i \in [K]$. This section constructs a single rule in $\cC$ that attains all of these bounds simultaneously, establishes simultaneous optimality and uniqueness of its conditional law, and compares it with e-BH and the randomized procedure U-eBH of \citet{Xu25}.

We begin by defining the candidate rule. Let $U \sim \operatorname{Uniform}(0, 1)$ be independent of $X$, and define
$$
D^\star(x, U) = \{i \in [K] : U \le G_i(x)\}.
$$
With a single uniform variable $U$ independent of $X$, write $\pi_i^\star(x) = \PP_U \{i \in D^\star(x, U)\}$ for the conditional inclusion probability. Because $U$ is uniform on $(0, 1)$, $\pi_i^\star(x) = 1 \wedge G_i(x)$. Attaining these bounds does not by itself place $D^\star$ in $\cC$; certification is verified in Theorem~\ref{thm-certification}.

The score $G_i(x)$ is monotone in $x_i$. If $x_i \ge x_j$, then $\min(x_i, x_\ell) \ge \min(x_j, x_\ell)$ for every $\ell$, so $G_i(x) \ge G_j(x)$. The rule therefore selects a prefix of the evidence ordering, and tied inputs are selected together. Zero inputs are never selected, as $G_i(x) = 0$ whenever $x_i = 0$. Because $D^\star$ selects a prefix, it is convenient to describe it through the ordered coordinates.

For $x_{[1]} \ge x_{[2]} \ge \dots \ge x_{[K]}$, define the ordered thresholds
$$
g_r = \frac{\alpha}{K} \left(r x_{[r]} + \sum_{j = r + 1}^K x_{[j]}\right), \:\: r \in [K],
$$
and the rejection count
$$
R^\star(x, U) = \max\left(\{r \in [K] : U \le g_r\} \cup \{0\}\right).
$$
The coordinate at rank $r$ satisfies $G_{[r]}(x) = g_r$: the terms with $j \le r$ contribute $x_{[r]}$ each, and the terms with $j > r$ contribute $x_{[j]}$. The thresholds satisfy $g_1 \ge g_2 \ge \dots \ge g_K \ge 0$; the nonnegativity is immediate from $x_{[j]} \ge 0$, and the monotonicity is verified as part of the proof of Theorem~\ref{thm-certification}. The set $D^\star(x, U)$ therefore consists of the coordinates at ranks up to $R^\star(x, U)$. We now show that $D^\star$ is conditionally certified.

\begin{theorem} \label{thm-certification}
Let $U \sim \operatorname{Uniform} (0, 1)$. For every $x \in [0, \infty)^K$ and every $i \in [K]$, $\EE_U [Z_i(D^\star(x, U))] \le x_i$. Consequently $D^\star \in \cC$, and Proposition~\ref{prop-fdr-null} gives $\fdr_P(D^\star) \le {\alpha K_0(P)} / {K}$ for every $P \in \cP$. 
\end{theorem}

\begin{proof}
Both $D^\star$ and the certificate $Z(D^\star)$ are equivariant under relabeling of the coordinates, so it suffices to prove the claim for a sorted input $x_1 \ge x_2 \ge \cdots \ge x_K$ and then extend by relabeling at the end. Fix $i \in [K]$ in this ordering.

Set $x_{K + 1} = 0$, $g_{K + 1} = 0$, and $h_r = 1 \wedge g_r$ for $r = 1, 2, \dots, K+1$. Using the definition of $g_r$, we write 
$$
g_r - g_{r + 1} = \frac{\alpha r}{K} (x_r - x_{r + 1}) \ge 0, \:\: r = 1, 2, \dots, K.
$$
Thus, $g_1 \ge g_2 \ge \dots \ge g_K \ge 0$ and $h_1 \ge h_2 \ge \dots \ge h_K \ge 0$. For $r = 1, 2, \dots, K$ define the intervals $I_r = (h_{r+1}, h_r]$. The intervals are disjoint with union $(0, h_1] \subseteq (0, 1]$, since $h_r \le 1$ for every $r$.

We claim that $R^\star(x, U) = r$ if and only if $U \in I_r$. Suppose first that $U \in I_r$, so that $h_{r + 1} < U \le h_r$. Since $h_r \le g_r$, we have $U \le g_r$ and therefore $R^\star \ge r$. To show that rank $r + 1$ fails, note that $U > h_{r + 1}$. If $g_{r + 1} > 1$, then $h_{r + 1} = 1$, which contradicts $U \le 1$; hence $g_{r + 1} \le 1$, so $h_{r + 1} = g_{r + 1}$ and $U > g_{r + 1}$. Since rank $r + 1$ fails the condition $U \le g_{r + 1}$, we have $R^\star < r + 1$, and combined with $R^\star \ge r$ this gives $R^\star = r$.

Conversely, suppose $R^\star = r$ with $r \ge 1$, so that $U \le g_r$ and $U > g_{r + 1}$. From $U > g_{r + 1}$ and $U \le 1$ we obtain $g_{r + 1} < 1$, hence $h_{r + 1} = g_{r + 1} < U$. For the upper endpoint, if $g_r \le 1$, then $h_r = g_r \ge U$, and if $g_r > 1$, then $h_r = 1 \ge U$; in either case $U \le h_r$. Together with $U > h_{r + 1}$ this gives $U \in I_r$. The case $R^\star = 0$ occurs if and only if $U > g_1$, which is equivalent to $U > h_1$ because $h_1 = 1 \wedge g_1$ and $U \le 1$. Therefore
$$
\PP_U(R^\star(x, U) = r) = h_r - h_{r + 1}, \:\: r = 1, 2, \dots, K,
$$
and $\PP_U(R^\star(x, U) = 0) = 1 - h_1$.

Returning to the fixed coordinate $i \in [K]$, we now compute its certificate expectation. For $r \ge i$, on the event $\{R^\star(x, U) = r\}$, the coordinate at rank $i$ belongs to $D^\star(x, U)$, so its certificate equals $Z_i(D^\star(x, U)) = K / (\alpha r)$. Consequently,
$$
\EE_U[Z_i(D^\star(x, U))]
= \sum_{r = i}^K \frac{K}{\alpha r} \PP_U(R^\star(x, U) = r)
= \sum_{r = i}^K \frac{K}{\alpha r} (h_r - h_{r+1}).
$$
Since $g_r \ge g_{r + 1}$ and clipping at one cannot increase a gap, $0 \le h_r - h_{r + 1} \le g_r - g_{r + 1}$. Substituting these bounds into the sum above gives
$$
\EE_U [Z_i(D^\star(x, U))]
\le \sum_{r = i}^K (x_r - x_{r+1})
= x_i - x_{K+1}
= x_i.
$$
For an arbitrary input $x$, let $\sigma$ be the permutation that sorts $x$ in decreasing order, so $x^\sigma_i = x_{\sigma(i)}=x_{[i]}$. The score $G_i(\cdot)$, the rule $D^\star$, and the certificate $Z_i(D^\star)$ are equivariant under relabeling, so the inequality proved for the sorted input transfers to $x$. Hence $\EE_U[Z_i(D^\star(x, U))] \le x_i$ for every $i \in [K]$. Therefore $D^\star \in \cC$, and the FDR bound is a direct consequence of Proposition~\ref{prop-fdr-null}.
\end{proof}

Combining the coordinatewise bound of Section~\ref{sec-certification-class} with Theorem~\ref{thm-certification} gives simultaneous optimality.
\begin{corollary} \label{cor-simultaneous}
For every $D \in \cC$, every $x \in [0, \infty)^K$, and every $i \in [K]$,
$$
\pi_i^D(x) \le \pi_i^\star(x) = 1 \wedge G_i(x).
$$
\end{corollary}

\begin{proof}
Section~\ref{sec-certification-class} gives $\pi_i^D(x) \le 1 \wedge G_i(x)$ for every $D \in \cC$. Theorem~\ref{thm-certification} gives $D^\star \in \cC$, and since $\pi_i^\star(x) = 1 \wedge G_i(x)$ by construction of $D^\star$, the claim follows.
\end{proof}

Corollary~\ref{cor-simultaneous} implies that every $D \in \cC$ attains at most the marginal inclusion probabilities of $D^\star(x, U)$ at this input. We now extend this to arbitrary weighted sums and show that any rule attaining the maximum with strictly positive weights has the same conditional law as $D^\star(x, U)$. We first record two consequences of equality in the marginal bounds. The following lemma records the equality conditions in the certificate bound.

\begin{lemma} \label{lem-equality-marginal}
Let $D$ satisfy $\EE_V [Z_i(D)] \le x_i$ for every $i \in [K]$. If $\pi_i^D(x) = G_i(x) \le 1$ for some $i \in [K]$, then $\EE_V [Z_i(D)] = x_i$ and $\EE_V [\min(Z_i(D), Z_j(D))] = \min(x_i, x_j)$ for every $j \in [K]$.
\end{lemma}

\begin{proof}
By the certificate identity and the certificate budget,
$$
\pi_i^D(x) = \frac{\alpha}{K} \sum_{j = 1}^K \EE_V [\min(Z_i(D), Z_j(D))] \le \frac{\alpha}{K} \sum_{j = 1}^K \min(x_i, x_j) = G_i(x).
$$
Since the two ends are equal, every term in the sum is tight, and the claim follows.
\end{proof}

\begin{lemma} \label{lem-equality-prefix}
Fix $x \in [0, \infty)^K$. Let $D$ satisfy $\EE_V [Z_i(D)] \le x_i$ and $\pi_i^D(x) = 1 \wedge G_i(x)$ for every $i \in [K]$. Then $D$ is almost surely a prefix of the evidence ordering, i.e, 
$$
j \in D \text{ and } x_i \ge x_j \implies i \in D \:\: \text{almost surely.}
$$
\end{lemma}

\begin{proof}
A coordinate with $G_i(x) > 1$ satisfies $\pi_i^D(x) = 1$ and is selected almost surely. A coordinate with $G_i(x) \le 1$ satisfies Lemma~\ref{lem-equality-marginal}. 

Fix $j \in [K]$ and suppose $j \in D$. If $G_j(x) > 1$, then $j$ is selected almost surely. For any $i$ with $x_i \ge x_j$, monotonicity of $G$ gives $G_i(x) \ge G_j(x) > 1$, so $i$ is also selected almost surely.

If $G_j(x) \le 1$, Lemma~\ref{lem-equality-marginal} applied to $j$ gives $\EE_V [Z_j(D)] = x_j$ and $\EE_V [\min(Z_j(D), Z_i(D))] = \min(x_j, x_i)$. For $i$ with $x_i \ge x_j$, we have $\min(x_j, x_i) = x_j = \EE_V [Z_j(D)]$, so
$\EE_V [\min(Z_j(D), Z_i(D))] = \EE_V [Z_j(D)]$. Since $\min(Z_j(D), Z_i(D)) \le Z_j(D)$ pointwise, the difference is nonnegative with zero expectation, hence zero almost surely. Therefore $Z_j(D) \le Z_i(D)$ almost surely. On the event $\{j \in D\}$, we have $Z_j(D) = K/(\alpha |D|) > 0$ whenever $|D| > 0$, so $Z_i(D) \ge Z_j(D) > 0$ and therefore $i \in D$. If $|D| = 0$, the case is vacuous.

Equal positive inputs satisfy $G_i(x) = G_j(x)$. If this common value is at least $1$, both are selected almost surely; otherwise, applying the argument above in both directions gives $i \in D \iff j \in D$ almost surely. Zero inputs have $G_i(x) = 0$ and are never selected. Therefore $D$ is almost surely a prefix of the evidence ordering.
\end{proof}

\begin{theorem} \label{thm-uniqueness}
Fix $x \in [0, \infty)^K$ and $w = (w_1, w_2, \dots, w_K) \in [0, \infty)^K$. For every random rejection set $D$ with $\EE_V [Z_i(D)] \le x_i$ for every $i \in [K]$,
\begin{equation} \label{eq-weighted-obj}
\sum_{i = 1}^K w_i \PP_V(i \in D) \le \sum_{i = 1}^K w_i \pi_i^\star(x).
\end{equation}
Moreover, if $w_i > 0$ for every $i$ and equality holds, then $D$ has the same conditional distribution as $D^\star$ at $x$.
\end{theorem}

\begin{proof}
By the oracle bound of Section~\ref{sec-certification-oracle} applied with $\cN = [K]$,
$$
\PP_V(i \in D) \le 1 \wedge G_i(x) = \pi_i^\star(x), 
$$
for every $i \in [K]$. Multiplying by $w_i \ge 0$ and summing over $i$ gives the stated inequality.

For the uniqueness claim, suppose $w_i > 0$ for every $i \in [K]$ and equality holds in \eqref{eq-weighted-obj}. Every term $w_i \PP_V (i \in D)$ is at most $w_i \pi_i^\star(x)$, which implies
$$
\PP_V(i \in D) = \pi_i^\star(x) = 1 \wedge G_i(x), 
$$
for every $i \in [K]$. By Lemma~\ref{lem-equality-prefix}, $D$ is almost surely a prefix of the evidence ordering. Write $\sigma_x(r)$ for the coordinate at rank $r$, so that $x_{\sigma_x(1)} \ge x_{\sigma_x(2)} \ge \cdots \ge x_{\sigma_x(K)}$, and set $h_r = 1 \wedge g_r$ for $r \in [K]$, as in the proof of Theorem~\ref{thm-certification}. For a prefix rule, $\sigma_x(r) \in D$ if and only if the rejection count $R = |D| \ge r$, so
$$
\PP (R \ge r) = \PP_V (\sigma_x(r) \in D) = 1 \wedge G_{\sigma_x(r)}(x) = h_r, \:\: r = 1, 2, \dots, K.
$$
$D$ is almost surely a prefix of the evidence ordering, i.e., $D$ is determined by $R$. The tail probabilities $\PP (R \ge r)$ determine the distribution of $R$ and therefore determine the distribution of $D$. The rule $D^\star$ is a prefix rule with the same tails, so $D$ and $D^\star$ have the same distribution at input $x$.
\end{proof}

\begin{corollary} \label{cor-unconditional}
For every $P \in \cP$, every weight vector $w(X) = (w_1(X), w_2(X), \dots, w_K(X))$ with $w_i(X) \ge 0$ measurable for every $i \in [K]$, and every $D \in \cC$,
$$
\EE_P \left[\sum_{i = 1}^K w_i(X) \bm 1 \{i \in D(X, V)\}\right] \le \EE_P \left[\sum_{i = 1}^K w_i(X) \bm 1\{i \in D^\star(X, U)\}\right].
$$
\end{corollary}

\begin{proof}
Fix $P \in \cP$ and $D \in \cC$. At each input $x$, Theorem~\ref{thm-uniqueness} applied with the weight vector $w(x)$ gives
$$
\sum_{i = 1}^K w_i(x) \pi_i^D(x) \le \sum_{i = 1}^K w_i(x) \pi_i^\star(x).
$$
Taking expectations over $X \sim P$ we get:
$$
\begin{aligned}
\EE_P \left[\sum_{i = 1}^K w_i(X) \bm 1 \{i \in D(X, V)\} \right]
=&\: \EE_P \left[\sum_{i = 1}^K w_i(X) \pi_i^D(X) \right] \\
\le&\: \EE_P \left[\sum_{i = 1}^K w_i(X) \pi_i^\star(X) \right] \\
=&\: \EE_P \left[\sum_{i = 1}^K w_i(X) \bm 1 \{i \in D^\star(X, U)\}\right].
\end{aligned}
$$
\end{proof}

Corollary~\ref{cor-unconditional} yields optimality for multiple objectives. With $w_i(X) = 1$, it gives optimality for the expected number of rejections. With $w_i(X) = \mathbf{1}\{i \notin \cN(P)\}$, it gives optimality for the expected number of true discoveries under any truth configuration; dividing by $K_1(P)$ extends this to average power when $K_1(P) > 0$. 

Corollary~\ref{cor-unconditional} establishes optimality of expected weighted counts. The comparison that follows is pathwise and holds for e-BH and U-eBH under a shared uniform seed.

Given $x$, draw $U \sim \operatorname{Uniform}(0, 1)$ and sort $x$ in decreasing order. The sample $D^\star(x, U)$ is obtained by computing the thresholds $g_r$ and returning the prefix of length $R^\star(x, U)$. The U-eBH rule of \citet{Xu25} applies e-BH to the scaled inputs $x_i / U$ for a uniform $U$ independent of $X$. Its rejection count is
$$
r_U(x, U) = \max\left(\left\{r \in [K] : U \le \frac{\alpha r}{K} x_{[r]}\right\} \cup \{0\}\right), 
$$
and the rejection set $D_U(x, U)$ consists of the coordinates at ranks up to $r_U(x, U)$. We use the same uniform seed $U$ for comparison.

\begin{proposition} \label{prop-containment}
For every $x \in [0, \infty)^K$ and every $U \in (0, 1)$,
$$
D_{\operatorname e}(x) \subseteq D_U(x, U) \subseteq D^\star(x, U).
$$
\end{proposition}

\begin{proof}
Let $r_{\operatorname e} = r_{\operatorname{e}}(x)$ denote the rejection rank of e-BH. If $r_e=0$ then the first inclusion above is trivial.  Otherwise, by definition of $r_{\operatorname e}$, $x_{[r_{\operatorname e}]} \ge K / (\alpha r_{\operatorname e})$, so $({\alpha r_{\operatorname e}} / {K}) x_{[r_{\operatorname e}]} \ge 1 > U$. Hence rank $r_{\operatorname e}$ satisfies the U-eBH condition, which implies $r_U \ge r_{\operatorname e}$ and therefore $D_{\operatorname e}(x) \subseteq D_U(x, U)$.

For the second inclusion, if rank $r$ satisfies $U \le ({\alpha r} / {K}) x_{[r]}$, then
$$
g_r = \frac{\alpha}{K} \left(r x_{[r]} + \sum_{j = r + 1}^K x_{[j]}\right) \ge \frac{\alpha r}{K} x_{[r]} \ge U. 
$$
This implies $R^\star \ge r_U$ and $D_U(x, U) \subseteq D^\star(x, U)$.
\end{proof}

Since $g_i \ge ({\alpha r} / {K}) x_{[r]}$ for every $r \ge i$, the inclusion probability of $D^\star$ at rank $i$ is at least that of $D_U$, i.e., 
$$
\PP_U \left(\sigma_x(i) \in D^\star(x, U) \right) - \PP_U \left(\sigma_x(i) \in D_U(x, U) \right) = (1 \wedge g_i) - (1 \wedge a_i(x)) \ge 0,
$$
where $a_i(x) = \max_{r \ge i} ({\alpha r} / {K}) x_{[r]}$. 

\section{Simulations} \label{sec-simulation}

This section reports a simulation study comparing the proposed rule with e-BH and with the competing randomized procedures. The goals are to measure the power improvement across signal strengths and dependence settings, and to verify that the attained FDR of every procedure remains below its theoretical bound.

\subsection{Simulation Setup}\label{sec-simulation-setup}

We consider a Gaussian mean problem with $K = 100$ hypotheses and target FDR $\alpha = 0.05$. In each replication an alternative set $\cA$ is drawn uniformly among the subsets of size $K_1 = 30$, independently of the noise, and the evidence is generated conditionally as 
$$
Z \mid \cA \sim N_K(\theta(\cA), \Sigma), \:\: \theta_i(\cA) = \mu \bm{1} \{i \in \cA\}.
$$
The hypotheses are $H_i : \theta_i \le 0$ against $\theta_i > 0$. True nulls have $\theta_i = 0$, at the boundary of the null hypothesis, and every alternative carries the same mean $\mu > 0$.

The covariance matrix $\Sigma$ takes one of two forms:
$$
\Sigma_{ij} =
\begin{cases}
\rho^{|i - j|}, & \text{Positive Toeplitz Dependence}, \\
1, & \text{Negative Equicorrelation}, \:\: i = j, \\
- \rho / (K - 1), & \text{Negative Equicorrelation}, \:\: i \ne j.
\end{cases}
$$
The parameter $\rho$ has different interpretations in the two cases: it is the adjacent-coordinate correlation in the Toeplitz family, and it controls the magnitude of the negative off-diagonal entries in the equicorrelation family. The signal strength $\mu$ takes equally spaced values on $[1, 4]$ in increments of $0.2$, and the dependence parameter $\rho$ takes equally spaced values on $[0,  0.9]$ in increments of $0.1$. Each setting is evaluated over $N = 5000$ Monte Carlo replications. The evidence supplied to the procedures is
$$
X_i = \exp \left(\lambda Z_i - \frac{\lambda^2}{2} \right), \:\: \lambda = \sqrt{2 \log \left(\frac{K}{2 \alpha} \right)}.
$$
The tuning parameter $\lambda$ is held fixed across the signal-strength grid, so the evidence construction is identical at every $\mu$.

Seven e-value procedures are compared:
\begin{itemize}
\item e-BH, the benchmark of Section~\ref{sec-setup-ebh},
\item R1-eBH, R2-eBH, and \Rboth{}, the three rounding procedures,
\item U-eBH and J-eBH, the shared-uniform and coordinate-specific randomizations,
\item $D^\star$, the procedure developed in Section~\ref{sec-optimality}.
\end{itemize}
The five randomized competitors are from \citet{Xu25}; their complete specifications are given in the appendix.

U-eBH and $D^\star$ share a single uniform $U$, which realizes the coupling $D_U(x,u) \subseteq D^\star(x, u)$ of Proposition~\ref{prop-containment} for every $x$ and $u$. The coordinatewise uniforms used by \Rone, \Rtwo, \Rboth, and \Jebh{} are independent across hypotheses, and the two stages of \Rboth{} use independent randomness. 

Let $X^{(b)}$ denote the evidence and $V_b$ the randomization in replication $b$, let $\cA_b$ be the alternative set in that replication, and let $\cN_b = [K] \setminus \cA_b$. The realized false discovery proportion and true-discovery proportion of a rule $D$ are
$$
F_b^D = \frac{|D(X^{(b)}, V_b) \cap \cN_b|}{|D(X^{(b)}, V_b)| \vee 1},
\:\:
T_b^D = \frac{|D(X^{(b)}, V_b) \cap \cA_b|}{K_1}.
$$
The reported FDR and power are the averages $\hat F^D = N^{- 1} \sum_{b = 1}^N F_b^D$ and $\hat T^D = N^{- 1} \sum_{b = 1}^N T_b^D$. Pointwise intervals use $\hat \xi \pm 1.96 s_\xi / \sqrt{N}$, where $\hat\xi$ and $s_\xi$ are the sample mean and standard deviation of the corresponding replication-level quantity and $s_\xi/\sqrt{N}$ is the Monte Carlo standard error. For $D^\star$ against a competitor $D$, the reported difference is the replication-level average of $\delta_b = T_b^{D^\star} - T_b^D$, with standard error computed from the $\delta_b$. Differences and their intervals are reported in percentage points.

\subsection{Simulation Results} \label{sec-simulation-results}

\begin{figure}[!ht]
\centering
\includegraphics[width=\textwidth]{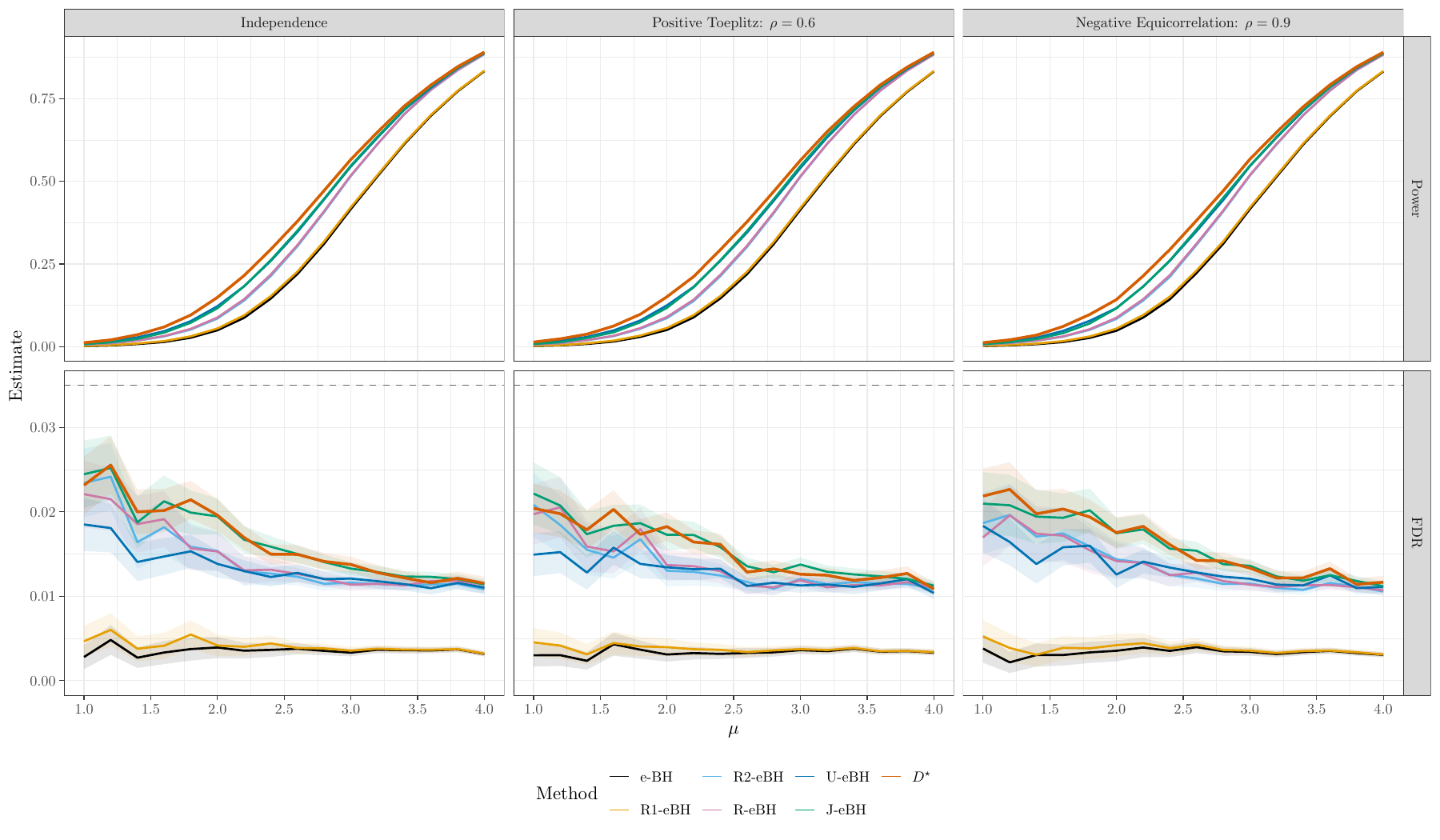}
\caption{Power (top row) and FDR (bottom row) against $\mu$ at independence ($\rho = 0$), positive Toeplitz dependence ($\rho = 0.6$), and negative equicorrelation ($\rho = 0.9$). Each panel compares the seven e-value procedures. Ribbons are pointwise $95\%$ Monte Carlo intervals. The dashed line in the lower panels marks the reference bound $\alpha K_0 / K = 0.035$.}
\label{fig-main-power-fdr}
\end{figure}

Figure~\ref{fig-main-power-fdr} reports power and FDR against $\mu$ at three representative dependence settings. Across the full grid, $D^\star$ attains the highest estimated power at every setting. The smallest gap over the grid is $0.03$ percentage points, attained against \Jebh{} at the extremes of the signal-strength grid, where every procedure is near zero power or near saturation. The largest is $16.4$ percentage points, attained against \ebh{} at intermediate signal strengths. The dominance of $D^\star$ over each competitor follows from the optimality theorem of Section~\ref{sec-optimality}. The simulation quantifies the magnitude of the gaps, their concentration in the $(\mu, \rho)$ grid, and the ordering among the competitors.

Power increases with $\mu$ for every procedure, and the separation between procedures is largest at intermediate signal strengths. \Uebh{} and \Jebh{} are the strongest competitors. At independence ($\rho = 0$) and $\mu = 2.4$, the power of $D^\star$ is $29.41\%$, compared with $25.95\%$ for \Uebh, $26.18\%$ for \Jebh, and $14.55\%$ for \ebh. The ordering is the same in the positive Toeplitz and negative equicorrelation panels of the figure.

\begin{figure}[!ht]
\centering
\includegraphics[width=\textwidth]{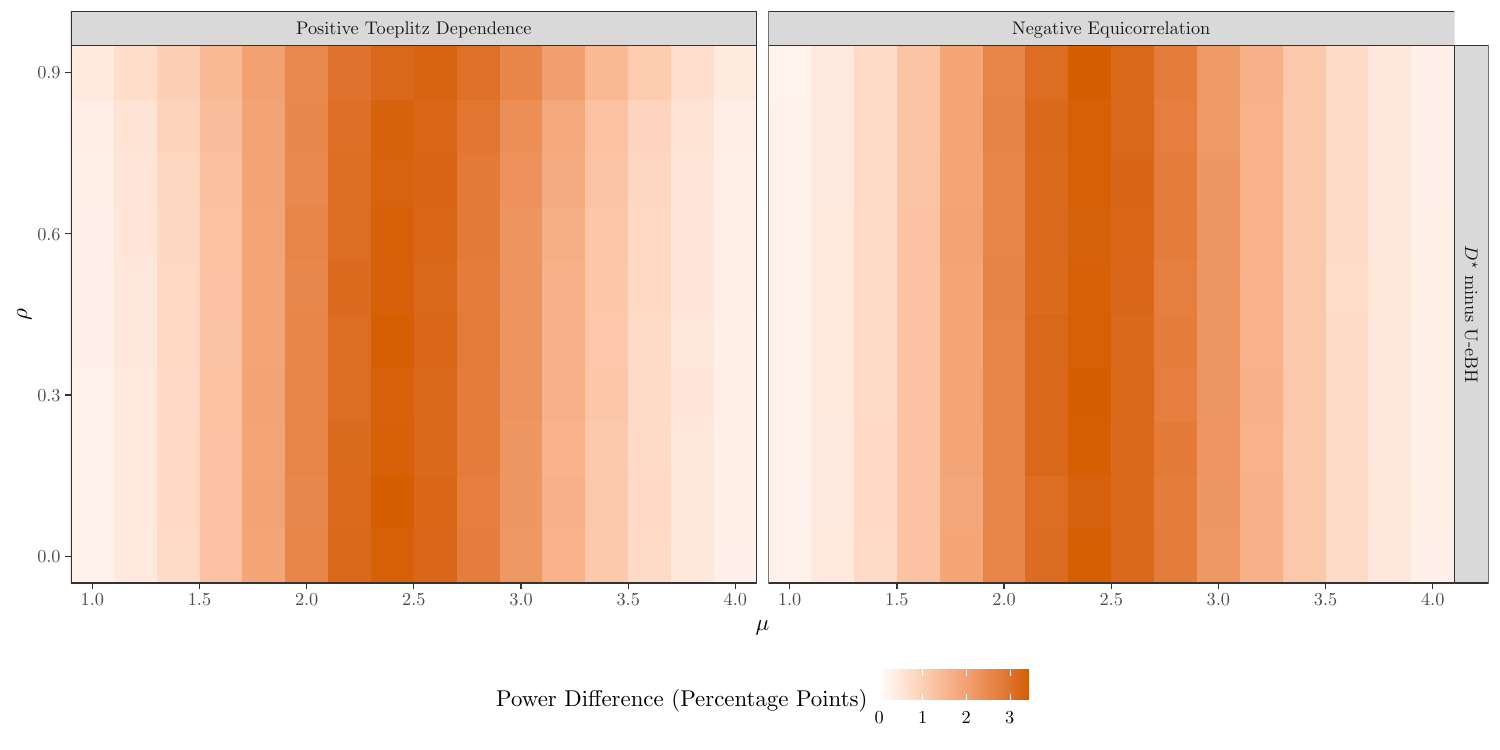}
\caption{Estimated power difference between $D^\star$ and \Uebh, over the $(\mu, \rho)$ grid.}
\label{fig-main-power-gain}
\end{figure}

Figure~\ref{fig-main-power-gain} reports the estimated difference in conditional expected power between $D^\star$ and \Uebh{} over the full grid, in percentage points. The gain is positive at every setting. It is concentrated at intermediate signal strengths, roughly $\mu \in [2.2, 2.8]$, and decays toward both endpoints: at $\mu = 1$ the gain is $0.27$ percentage points and at $\mu = 4$ it is $0.30$ points, against a peak of $3.43$ points. The peak is attained at $(\mu, \rho) = (2.4, 0.1)$ in the positive family and at $(\mu, \rho) = (2.4, 0.9)$ in the negative family. 

\begin{figure}[!ht]
\centering
\begin{subfigure}{0.49 \textwidth}
\includegraphics[width= \textwidth]{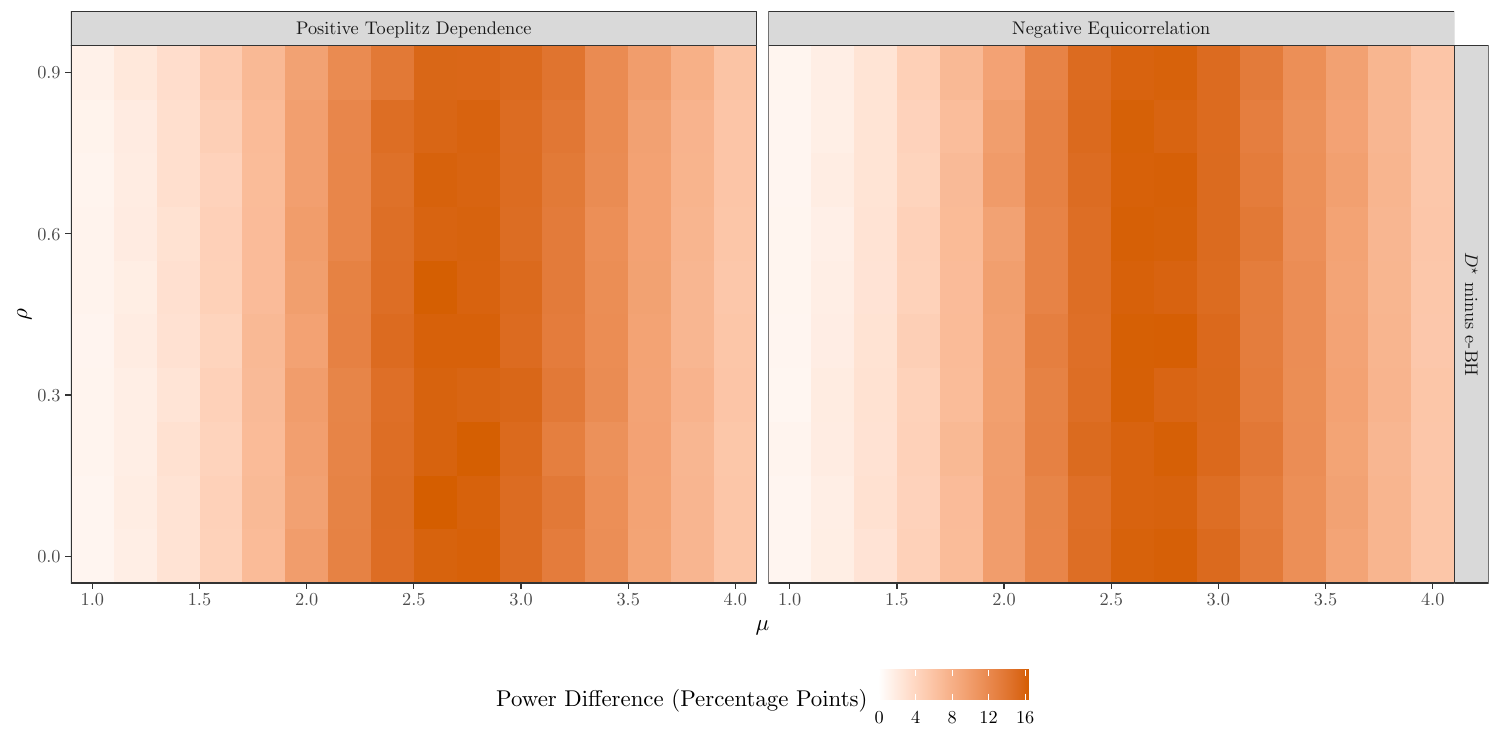}
\caption{$D^\star$ versus \ebh}
\end{subfigure}\hfill
\begin{subfigure}{0.49 \textwidth}
\includegraphics[width=\textwidth]{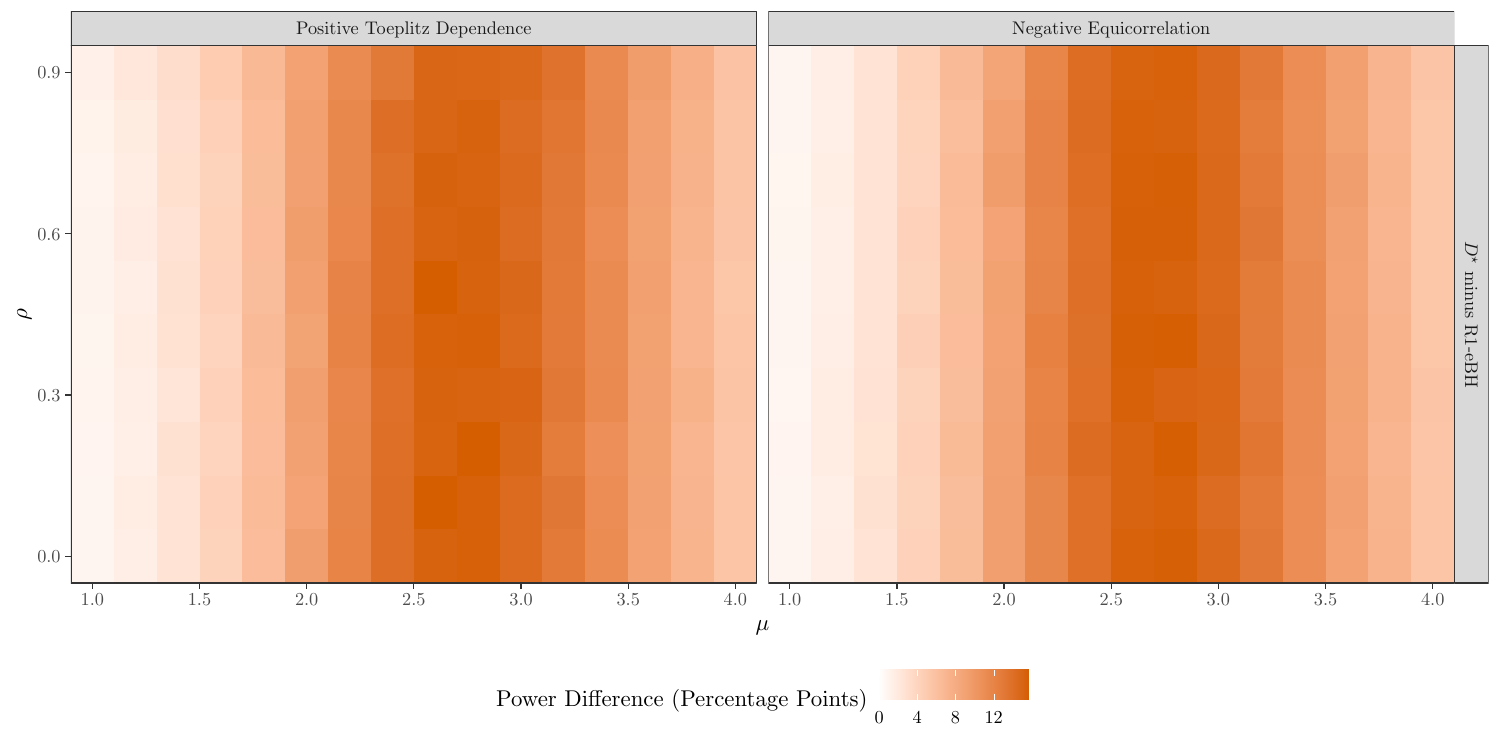}
\caption{$D^\star$ versus \Rone}
\end{subfigure}

\vspace{0.4cm}
\begin{subfigure}{0.49 \textwidth}
\includegraphics[width=\textwidth]{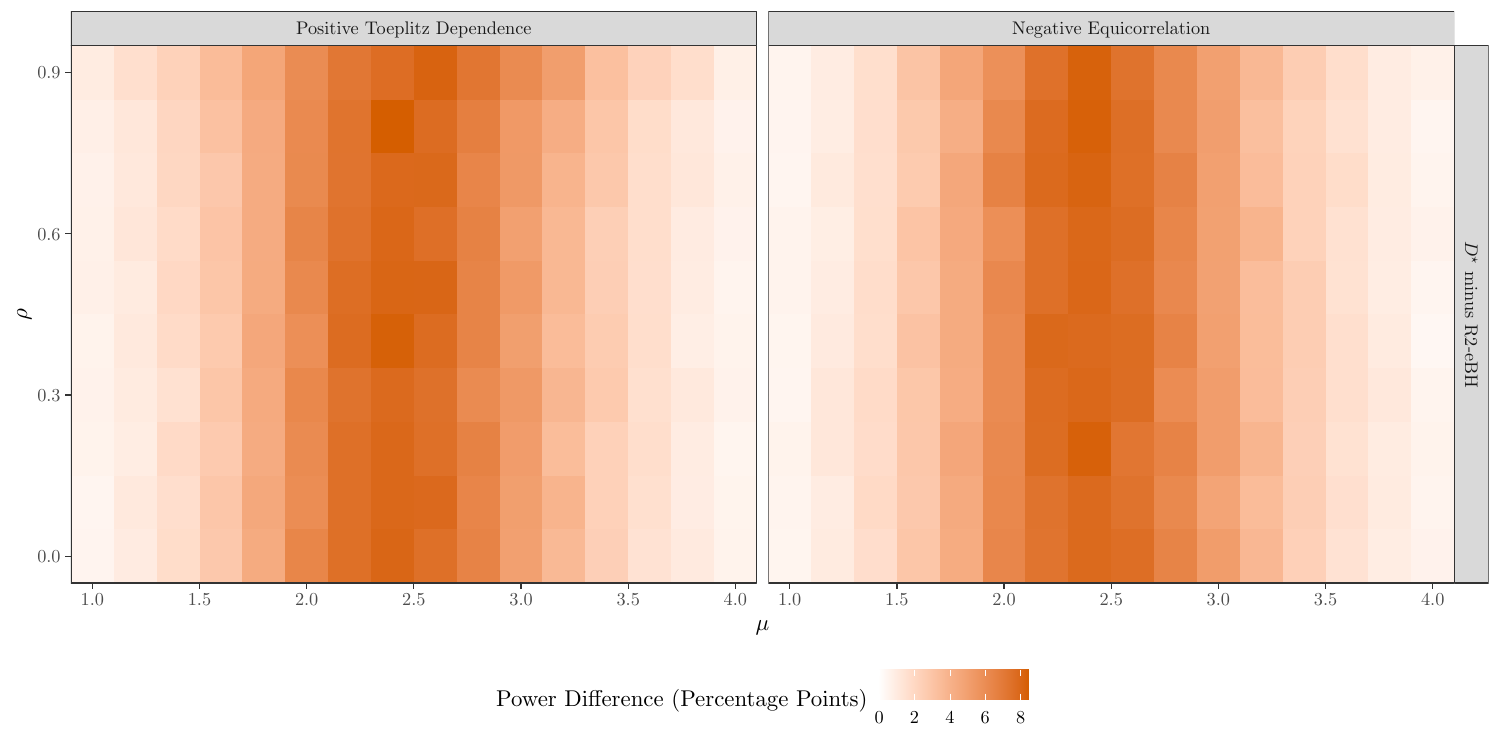}
\caption{$D^\star$ versus \Rtwo}
\end{subfigure}\hfill
\begin{subfigure}{0.49 \textwidth}
\includegraphics[width=\textwidth]{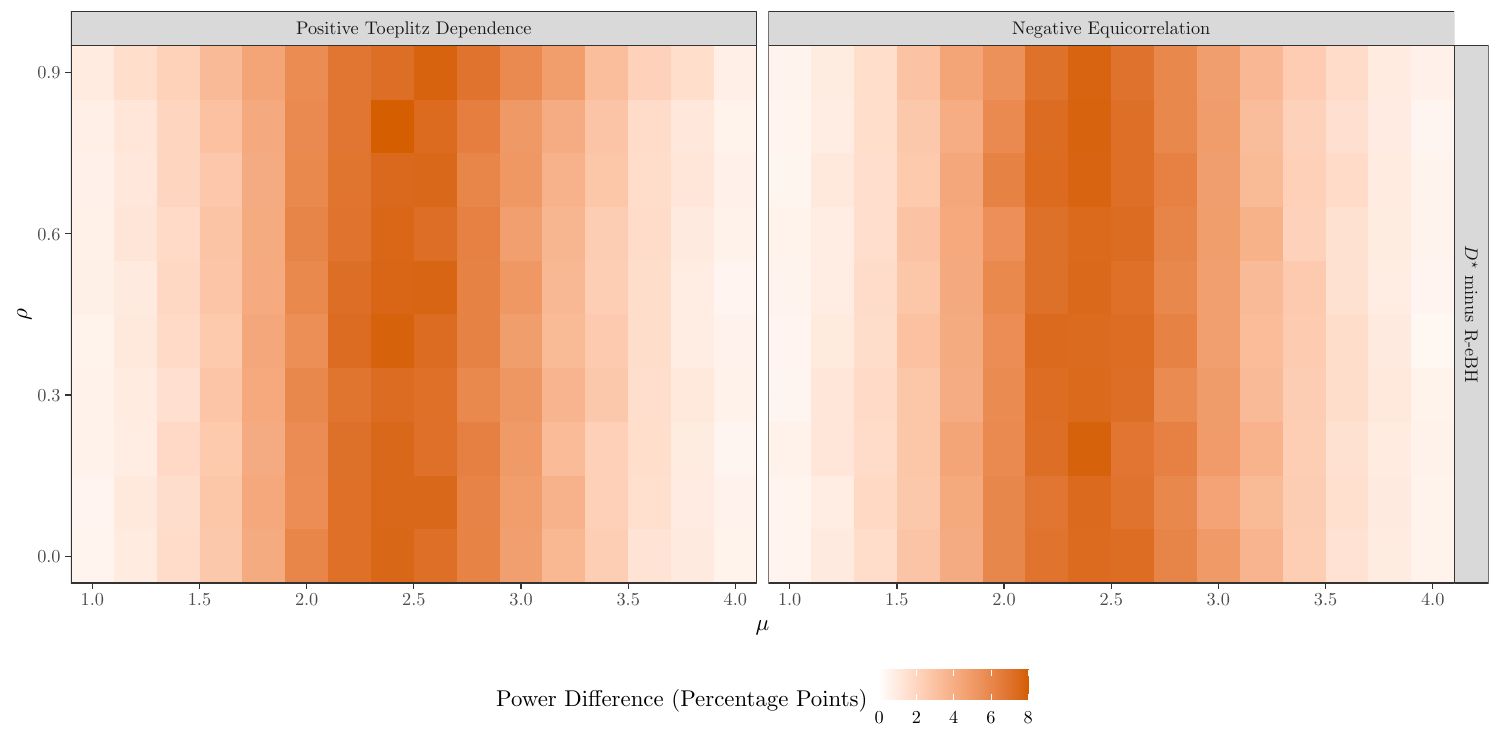}
\caption{$D^\star$ versus \Rboth}
\end{subfigure}

\vspace{0.4cm}
\begin{subfigure}{0.49 \textwidth}
\includegraphics[width=\textwidth]{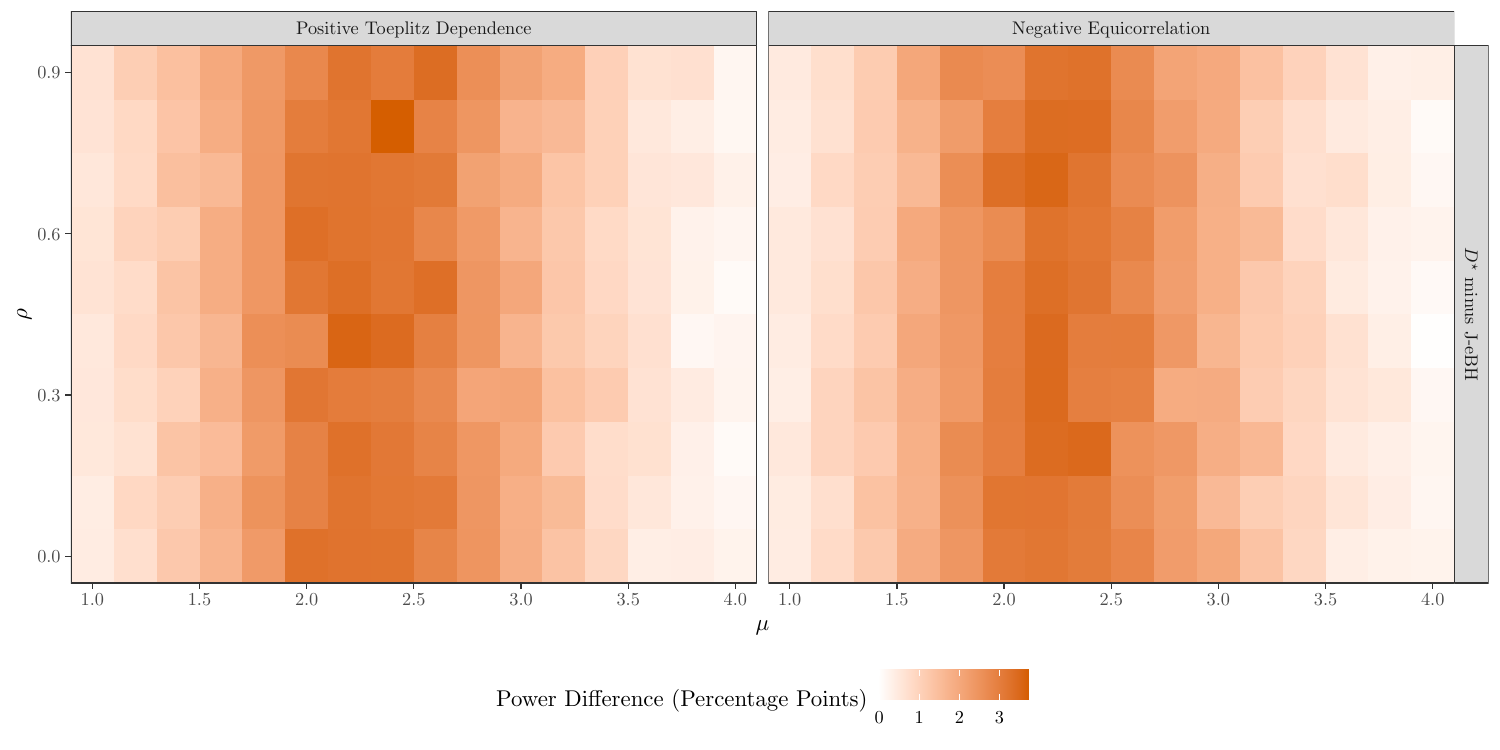}
\caption{$D^\star$ versus \Jebh}
\end{subfigure}
\caption{Paired power differences of $D^\star$ vs. each competitor, over the $(\mu, \rho)$ grid. Color intensities are not comparable across panels.}
\label{fig-supp-power-heatmaps}
\end{figure}

Figure~\ref{fig-supp-power-heatmaps} reports the corresponding heatmaps for $D^\star$ against the remaining e-value competitors. The maxima over the grid are $16.4$ percentage points over \ebh, $15.7$ over \Rone, $8.5$ over \Rtwo, $8.0$ over \Rboth, and $3.7$ over \Jebh. The medians of these gains over the positive Toeplitz settings are $9.5$, $9.2$, $3.5$, $3.4$, and $1.6$ percentage points for the same five competitors. Each median is at least $40\%$ of the corresponding maximum, so the large gains are not confined to the settings where the maxima are achieved. Each heatmap shows the same pattern: a broad intermediate-signal region of improvement, with the gains decaying toward both endpoints. The intermediate-signal improvement persists in both covariance families and in both cases the heatmaps vary relatively little with $\rho$.

\begin{figure}[!ht]
\centering
\begin{subfigure}{0.49 \textwidth}
\includegraphics[width=\textwidth]{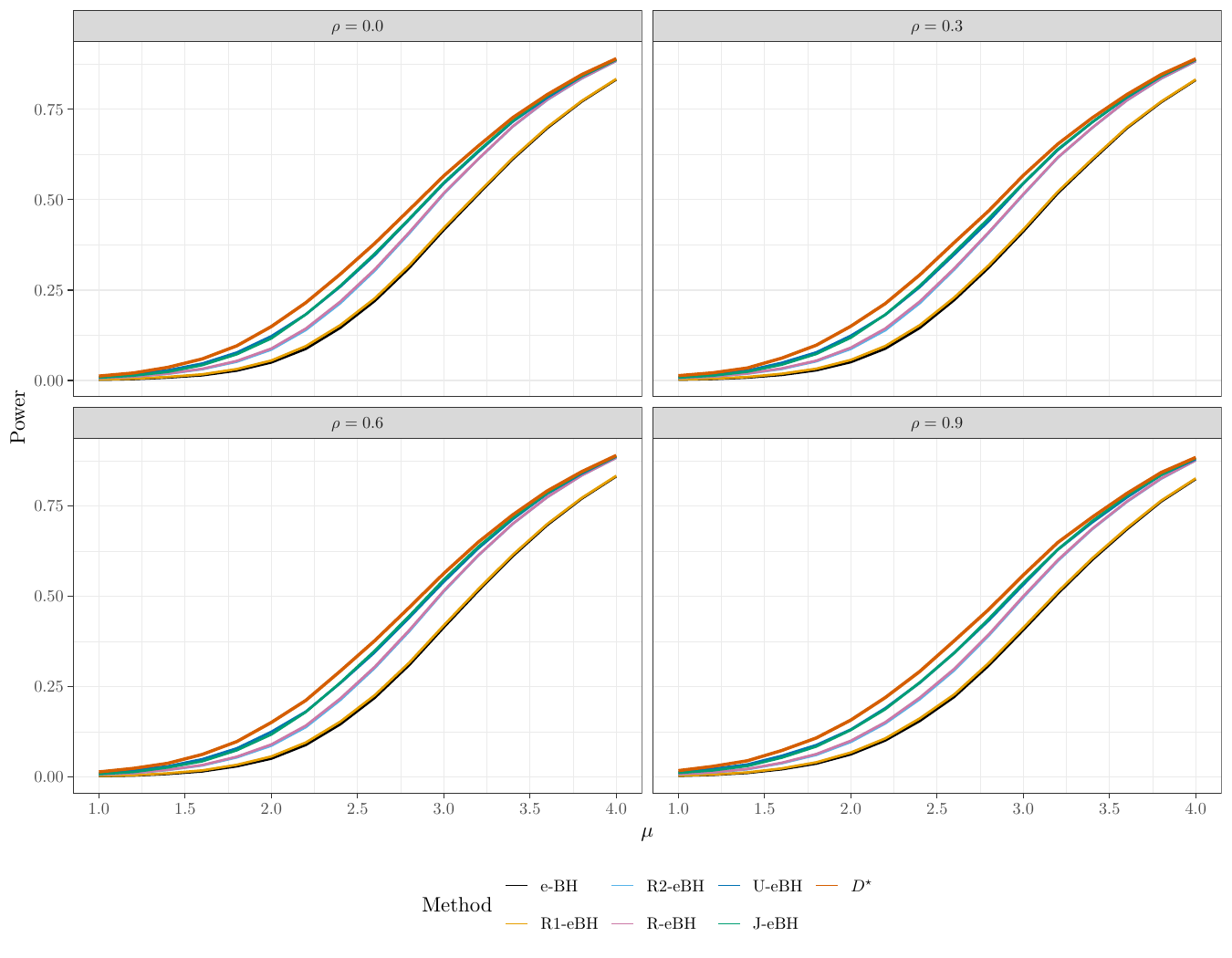}
\caption{Power, Positive Toeplitz}
\end{subfigure}\hfill
\begin{subfigure}{0.49 \textwidth}
\includegraphics[width=\textwidth]{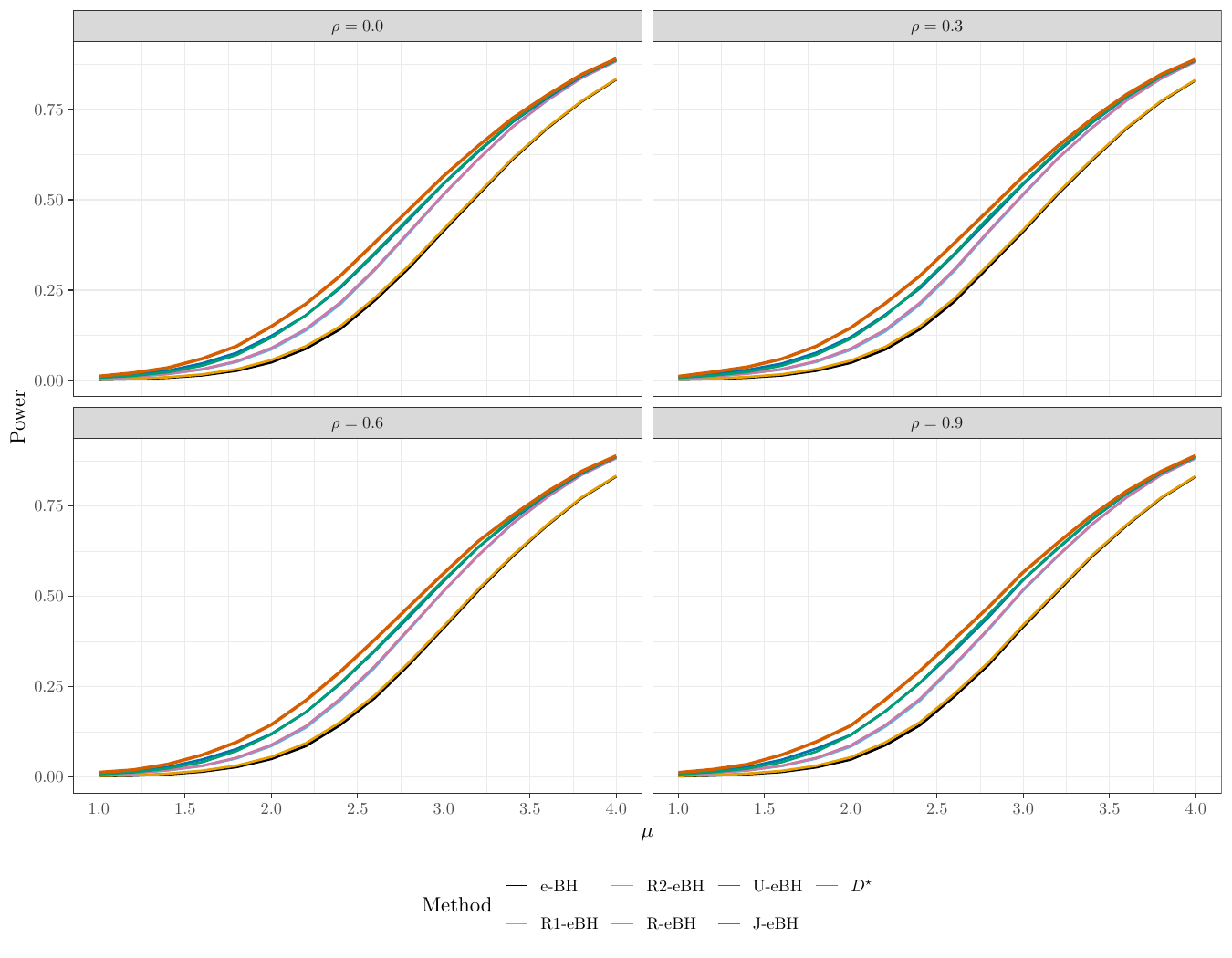}
\caption{Power, Negative equicorrelation}
\end{subfigure}

\vspace{0.4cm}
\begin{subfigure}{0.49 \textwidth}
\includegraphics[width=\textwidth]{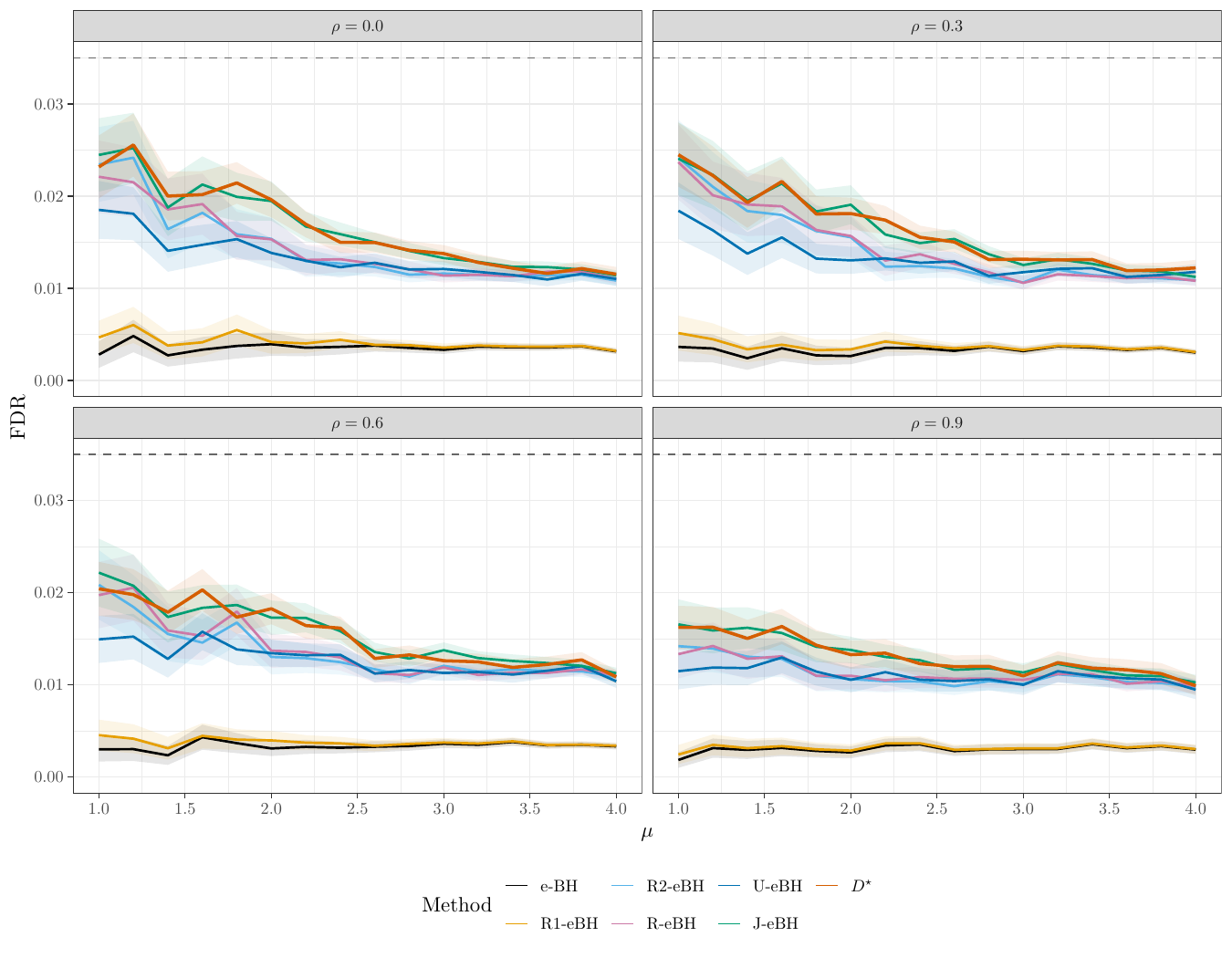}
\caption{FDR, Positive Toeplitz}
\end{subfigure}\hfill
\begin{subfigure}{0.49 \textwidth}
\includegraphics[width=\textwidth]{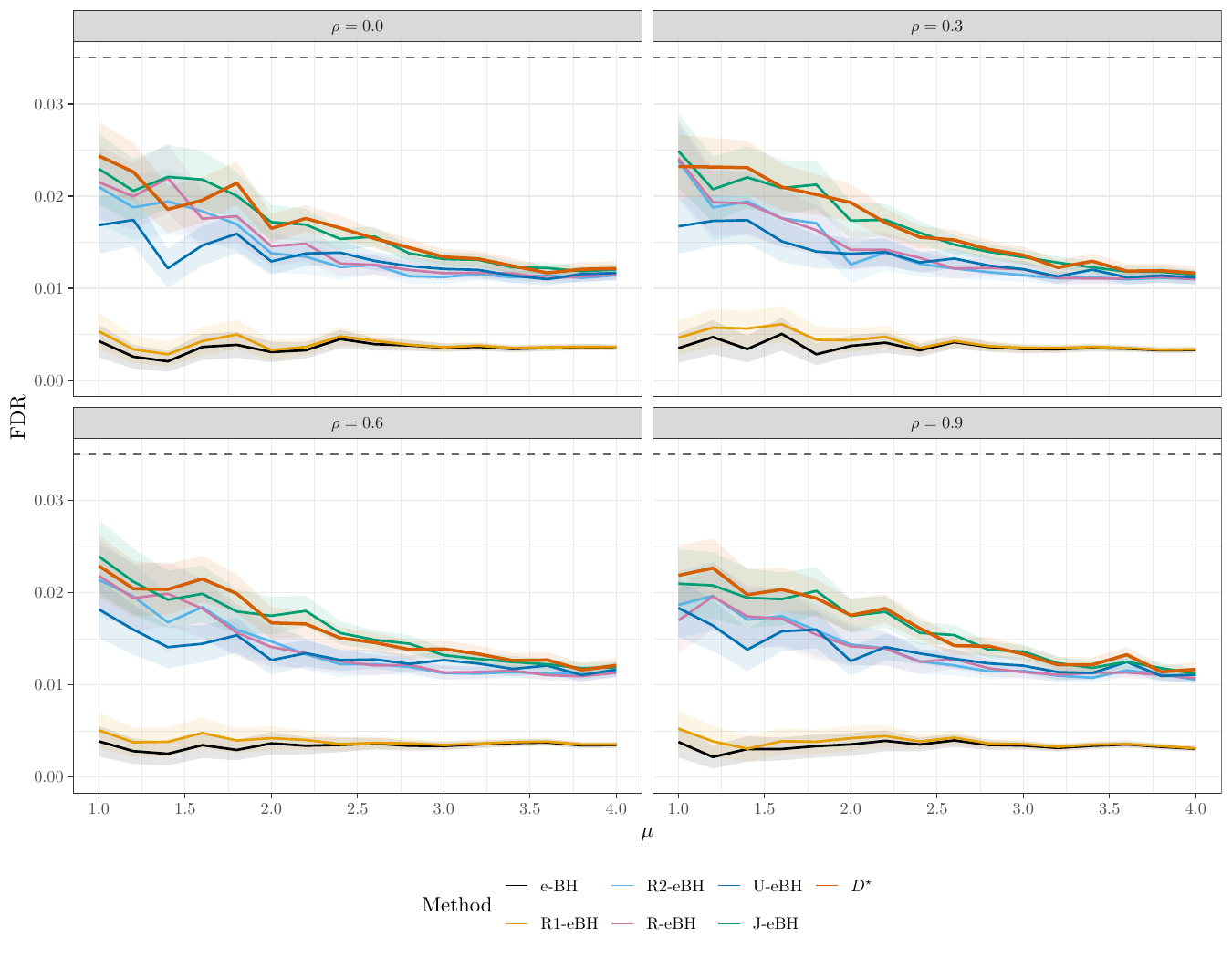}
\caption{FDR, Negative equicorrelation}
\end{subfigure}
\caption{Power and FDR curves for all seven e-value procedures at $\rho \in \{0, 0.3, 0.6, 0.9\}$ in each dependence family. Ribbons are pointwise $95\%$ Monte Carlo intervals; the dashed line in the FDR panels marks $\alpha K_0 / K = 0.035$.}
\label{fig-supp-seven-method-curves}
\end{figure}

Figure~\ref{fig-supp-seven-method-curves} reports the complete power and FDR curves for all seven e-value procedures at $\rho \in \{0, 0.3, 0.6, 0.9\}$ in each dependence family. The power ordering established by Figure~\ref{fig-main-power-fdr} persists across all eight panels. The FDR of every procedure lies below $0.035$ throughout. 

\section{Discussion} \label{sec-discussion}

This paper characterizes the simultaneous inclusion-probability optimization for e-BH under conditional certification and solves it. The certificate identity yields a coordinatewise upper bound, a single rule with one shared uniform variable attains every bound at once, and the equality conditions force any such rule to reject a prefix of the evidence ordering and fix its conditional rejection-set distribution. Simultaneous attainment implies optimality for every nonnegative weighted sum of inclusion probabilities, including expected true discoveries under any configuration of true and false nulls. The rule requires no additional tuning parameter beyond the target level, and its rejection set contains those of e-BH and U-eBH pathwise under the same uniform draw.

The optimization is relative to the certified class, which imposes budgets coordinate-by-coordinate, at each input on every coordinate, and membership is relative to the chosen evidence representation. Frameworks based on compound e-values \citep{Ignatiadis24} or e-closure \citep{Xu25b} can permit procedures that need not satisfy these budgets, and the results do not constrain them. The result establishes optimality within the certified class, not admissibility among all FDR-controlling rules. The complete classes of \citet{Sun26} address a separate question. The BY construction in the appendix optimizes over the calibrated evidence, so its optimality is relative to that representation rather than to all procedures using the original p-values.

Certification can be weakened to hold $P$-almost surely. For each fixed distribution, budgets on the true nulls suffice for FDR control, and budgets on all coordinates preserve the inclusion bounds and power comparisons almost surely. The pointwise form is retained because it defines the optimization at every input without reference to a model.

Two directions remain. The certificate identity relies on a common certificate value on the rejection set, zero elsewhere, and a rejection count entering through its reciprocal. Whether other error criteria admit analogous identities and budgets is open. The false coverage rate is a candidate, but an optimality characterization for it would differ from the randomized improvement already developed for that criterion \citep{Xu25}. The results here optimize the decision given the evidence and an independent randomization seed. Optimizing evidence construction jointly with the decision under a specified model is a separate problem that this analysis does not address.

\newpage{}


\newpage{}

\appendix

\singlespacing

\setstretch{2}

\section{Appendix} \label{sec-appendix}

This appendix collects the definitions of the randomized e-BH procedures of \citet{Xu25}, together with the certificate bounds that place them in the class $\cC$ of conditionally certified rules. The bounds below restate that guarantee in the certificate language of Section~\ref{sec-certification}, which the optimality analysis requires at every coordinate rather than only in aggregate over the nulls. In each case the bound can be derived from a randomized evidence vector whose coordinatewise conditional expectation is the original evidence, together with the self-consistency of e-BH applied to that vector, which is the step underlying the FDR proof of \citet{Wang22}. Section~\ref{app-stochastic-rounding} and Section~\ref{app-uebh-jebh} collect the definitions, Section~\ref{app-certified-procedures} proves the certificate bounds, and Section~\ref{app-by-certified} extends the construction to BY-calibrated p-values, establishing a pathwise nesting among the ordinary, shared-uniform, and optimal rejection sets.

\subsection{Stochastic Rounding and Related Procedures} \label{app-stochastic-rounding}

Let $\cG \subseteq [0, \infty]$ be a closed set. For $x \in [0, \infty)$ let $g_L(x) = \sup\{g \in \cG : g \le x\}$ and $g_U(x) = \inf\{g \in \cG : g \ge x\}$, both of which lie in $\cG$ when $x$ is within the range of $\cG$. Stochastic rounding of $x$ onto $\cG$ is
$$
T_{\cG}(x, U) =
\begin{cases}
x, & x \in \cG,\ x < \inf \cG, \: x > \sup \cG, \:\text{or}\: g_U(x) = \infty, \\
g_U(x), & g_L(x) < x < g_U(x) < \infty \:\text{and}\: U \le \frac{x - g_L(x)}{g_U(x) - g_L(x)}, \\
g_L(x), & g_L(x) < x < g_U(x) < \infty \:\text{and}\: U > \frac{x - g_L(x)}{g_U(x) - g_L(x)},
\end{cases}
$$
with $U \sim \operatorname{Uniform}(0, 1)$. Elements of $\cG$ and values outside the range of $\cG$ are left unchanged. For $x$ strictly between adjacent grid points $a < b$, the probability of rounding to $b$ is $(x - a) / (b - a)$, so that
$$
\EE_U [T_{\cG}(x, U)] = g_L(x) \cdot \frac{g_U(x) - x}{g_U(x) - g_L(x)} + g_U(x) \cdot \frac{x - g_L(x)}{g_U(x) - g_L(x)} = x.
$$
Applying $T_{\cG}$ to an e-value $X$ with $U$ independent of $X$ therefore produces an e-value $T_{\cG}(X,U)$ with the same expectation. 

The grid used by e-BH is $\cK = \{0\} \cup \{q_r : r \in [K]\} \cup \{\infty\}$ where $q_r = K / (\alpha r)$ are the thresholds at which e-BH can reject. A coordinate already at a threshold is unchanged by rounding. The \Rone{} procedure applies e-BH to the rounded evidence on $\cK$, with independent uniform variables $U_1, U_2, \dots, U_K$: 
$$
D_1(x, U) = D_{\operatorname e}\left(T_{\cK} (x_1, U_1), T_{\cK} (x_2, U_2), \dots, T_{\cK} (x_K, U_K) \right).
$$

The \Rtwo{} procedure replaces the fixed grid by a data-dependent threshold. Write $r_{\operatorname e}(x) = |D_{\operatorname e} (x)|$ and set 
$$
\hat \alpha^*(x) = \alpha(r_{\operatorname e} (x) + 1) / K, \:\: \cG_\alpha = \{0, \left(\hat \alpha^* (x)\right)^{- 1}, \infty\}.
$$
The procedure applies stochastic rounding to each coordinate onto the grid $\cG_\alpha$ and then applies e-BH to the rounded vector:
$$
D_2(x, U) = D_{\operatorname e} \left(T_{\cG_\alpha}(x_1, U_1), T_{\cG_\alpha}(x_2, U_2), \dots, T_{\cG_\alpha}(x_K, U_K) \right).
$$

The \Rboth{} procedure composes the two roundings. Let $U'$ and $U''$ be independent uniform vectors, set $W = (T_{\cK} (x_1,U_1'), T_{\cK} (x_2,U_2') \dots, T_{\cK}(x_K, U_K'))$, and define
$$
D_{\operatorname{both}}(x, (U', U'')) = D_2 (W, U'').
$$

\subsection{Shared and Coordinate-Specific Randomization} \label{app-uebh-jebh}

The \Uebh{} procedure applies e-BH to the scaled evidence $x_i / U$ with a single uniform variable $U \sim \operatorname{Uniform} (0, 1)$ shared across all coordinates:
$$
D_U(x, U) = D_{\operatorname e} \left(x_1 / U, x_2 / U, \dots, x_K / U \right).
$$
Writing $x_{[1]} \ge x_{[2]} \ge \dots \ge x_{[K]}$ for the decreasing order of $x$, its rejection count is
$$
r_U(x, U) = \max \left(\{r \in [K] : U \le x_{[r]} / q_r\} \cup \{0\}\right),
$$
and $D_U(x, U)$ consists of the coordinates at ranks up to $r_U(x, U)$.

The \Jebh{} procedure uses independent uniform variables $U_1, U_2, \dots, U_K$ and applies the BH procedure at level $\alpha$ to the ratios $U_i / x_i$, with $U_i / x_i = \infty$ when $x_i = 0$. Its rejection count is
$$
k_J(x, U) = \max \left(\left\{ k \in [K] : \sum_{i \in [K]} \bm{1} \{U_i \le x_i / q_k \} \ge k \right\} \cup \{0\}\right),
$$
and $D_J(x, U) = \{i \in [K] : U_i \le x_i / q_{k_J}\}$, with $D_J(x, U) = \emptyset$ when $k_J = 0$.

\subsection{Conditional Certification of Randomized e-BH Procedures} \label{app-certified-procedures}

Recall the certificate $Z_i$ of Section~\ref{sec-certification-fdr}. This subsection establishes the certificate bound for the three rounding procedures of Section~\ref{app-stochastic-rounding} and for the two uniform-based procedures of Section~\ref{app-uebh-jebh}.

\begin{proposition} \label{prop-rounding-certified}
The procedures \Rone{}, \Rtwo{}, and \Rboth{} satisfy $\EE_U [Z_i(D(x, U))] \le x_i$ for every $x \in [0, \infty)^K$ and every $i \in [K]$.
\end{proposition}

\begin{proof}
Each of the three procedures has the form $D(x, U) = D_{\operatorname e} (W(x, U))$ for a jointly randomized nonnegative vector $W$, and in each case the definition of stochastic rounding gives $\EE_U [W_i(x, U)] = x_i$. Since e-BH is self-consistent with respect to its input, $D_{\operatorname e} (W)$ is self-consistent with respect to $W$: for every realization and every $i \in D_{\operatorname e} (W)$ we have $W_i \ge K / (\alpha |D_{\operatorname e} (W)|)$. The certificate $Z_i(D)$ is the coordinatewise smallest nonnegative vector supporting $D$, so $Z_i (D_{\operatorname e} (W)) \le W_i$ pointwise. Taking expectations and applying $\EE_U [W_i] = x_i$ gives $\EE_U [Z_i(D(x, U))] \le x_i$.

It remains to identify $W$ in the three cases. For \Rone{}, 
$$
W = (T_{\cK} (x_1, U_1), T_{\cK} (x_2,U_2), \dots, T_{\cK} (x_K, U_K))
$$
with independent $U_i$, and $\EE_U [W_i] = x_i$ by the mean-preservation identity of stochastic rounding on $\cK$. For \Rtwo{}, 
$$
W = (T_{\cG_\alpha} (x_1, U_1), T_{\cG_\alpha} (x_2, U_2), \dots, T_{\cG_\alpha} (x_K, U_K))
$$
and the same identity on $\cG_\alpha$ gives $\EE_U [W_i] = x_i$. For \Rboth{}, let 
$$
\begin{aligned}
W' =&\: (T_{\cK} (x_1, U_1'), T_{\cK} (x_2, U_2'), \dots, T_{\cK} (x_K, U_K')), \\
W'' =&\: (T_{\cG_\alpha (W')} (W_1', U_1''), T_{\cG_\alpha (W')} (W_2', U_2''), \dots, T_{\cG_\alpha (W')} (W_K', U_K''))
\end{aligned}
$$
with $U''$ independent of $U'$. Conditional on $W'$, the mean-preservation identity on $\cG_\alpha (W')$ gives $\EE_{U''}[W_i'' \mid W'] = W_i'$, so $\EE [W_i''] = \EE [W_i'] = x_i$.
\end{proof}

\begin{proposition}\label{prop-uniform-certified}
The procedures \Uebh{} and \Jebh{} satisfy $\EE_U [Z_i(D(x, U))] \le x_i$ for every $x \in [0, \infty)^K$ and every $i \in [K]$.
\end{proposition}

\begin{proof}
For \Uebh{}, fix a rank $r \in [K]$, write $i_r$ for the coordinate at rank $r$ so that $x_{i_r} = x_{[r]}$, and write $Z_{[r]} (A) = Z_{i_r} (A)$ for the certificate of that coordinate. Choose $\ell \in \arg\max_{s \ge r} x_{[s]}/q_s$. Rank $r$ is rejected when some rank $s \ge r$ is feasible, which is equivalent to
$$
r_U \ge r \:\: \Longleftrightarrow \:\: U \le \frac{x_{[\ell]}}{q_\ell}.
$$
On this event, rank $\ell$ is feasible, so $r_U \ge \ell$. Therefore
$$
Z_{[r]} (D_U) = \frac{K}{\alpha (r_U \vee 1)} \bm{1} \{r_U \ge r\} \le \frac{K}{\alpha \ell} \bm{1} \left\{U \le \frac{x_{[\ell]}}{q_\ell}\right\}.
$$
Taking expectations and using the uniformity of $U$, we have 
$$
\EE_U [Z_{[r]} (D_U)] \le q_\ell \min \left(1, \frac{x_{[\ell]}}{q_\ell} \right) = \min(q_\ell, x_{[\ell]}) \le x_{[\ell]} \le x_{[r]}. 
$$

For \Jebh{}, write $k_J = |D_J|$.  The $k_J=0$ case is trivial, so assume that $k_J\ge 1$. For each coordinate $i \in [K]$, define
$$
L_i (x, U_{- i}) = \max \left(\left\{r \in [K] : 1 + \sum_{j \neq i} \bm{1} \{U_j \le x_j / q_r\} \ge r \right\} \cup \{0\} \right).
$$
This is the step-up count when coordinate $i$ is forced to be selected. By definition $L_i \ge 1$. The sum runs over $j \neq i$, so $L_i$ is a function of $U_{- i} = \left(U_1, U_2, \dots, U_{i - 1}, U_{i + 1}, \dots, U_K \right)$ alone. Since $k_J$ is the rejection count, at threshold $q_{k_J}$ the full sum satisfies $\sum_{j \in [K]} \bm{1} \{U_j \le x_j / q_{k_J}\} \ge k_J$. This implies $1 + \sum_{j \neq i} \bm{1} \{U_j \le x_j / q_{k_J}\} \ge k_J$. Rank $k_J$ is therefore feasible in the definition of $L_i$, and $L_i \ge k_J$.

We show that $i \in D_J$ if and only if $U_i \le x_i / q_{L_i}$, and that $k_J = L_i$ on $\{i \in D_J\}$. If $U_i \le x_i / q_{L_i}$, the definition of $L_i$ gives $1 + \sum_{j \neq i} \bm{1} \{U_j \le x_j / q_{L_i}\} \ge L_i$. Since $i$'s contribution equals one, the full sum satisfies $\sum_{j \in [K]} \bm{1} \{U_j \le x_j / q_{L_i}\} \ge L_i$, so at least $L_i$ coordinates satisfy their threshold at rank $L_i$ and hence $k_J \ge L_i$. Combined with $k_J \le L_i$ this gives $k_J = L_i$. The hypothesis is also the feasibility condition at rank $k_J$, since $q_{L_i} = q_{k_J}$, so $i \in D_J$. Conversely, if $i \in D_J$, then $U_i \le x_i / q_{k_J}$. Since $k_J \le L_i$ and $q_s = K / (\alpha s)$ is decreasing in $s$, $q_{k_J} \ge q_{L_i}$, so $x_i / q_{k_J} \le x_i / q_{L_i}$ and hence $U_i \le x_i / q_{L_i}$. Therefore $i \in D_J \iff U_i \le x_i / q_{L_i}$, and hence $k_J = L_i$ whenever $i \in D_J$.

On $\{i \in D_J\}$ the certificate is $K / (\alpha k_J)$ and zero elsewhere. Since $k_J = L_i$ on that event and $i \in D_J \iff U_i \le x_i / q_{L_i}$, the certificate is $q_{L_i}$ on $\{U_i \le x_i / q_{L_i}\}$ and zero elsewhere, so
$$
Z_i(D_J) = q_{L_i} \bm{1} \left\{U_i \le \frac{x_i}{q_{L_i}}\right\}.
$$
The quantity $L_i$ is a function of $U_{- i}$ at fixed $x$, so
$$
\EE_{U_i} [Z_i(D_J) \mid U_{- i}] = q_{L_i} \PP \left(U_i \le \left. \frac{x_i}{q_{L_i}} \right| U_{- i}\right) = q_{L_i} \min \left(1, \frac{x_i}{q_{L_i}}\right) \le x_i.
$$
Taking expectations over $U_{-i}$ gives $\EE_U [Z_i(D_J)] \le x_i$. 
\end{proof}

\subsection{Certified BY-Calibrated Procedures} \label{app-by-certified}

For $r \in [K]$, write $H_K = \sum_{r = 1}^K r^{-1}$ and $t_r = \alpha r / (K H_K)$, with $t_0 = 0$. The Benjamini-Yekutieli calibration function is
$$
f_{\operatorname{BY}} (p) =
\begin{cases}
{K} / (\alpha r), & t_{r - 1} < p \le t_r, \:\:r \in [K], \\
0, & p > t_K.
\end{cases}
$$
$f_{\operatorname{BY}}$ is nonnegative and nonincreasing, and
$$
\int_0^1 f_{\operatorname{BY}} (p) dp = \sum_{r = 1}^K \frac{K}{\alpha r} \cdot \frac{\alpha}{K H_K} = \frac{1}{H_K} \sum_{r = 1}^K \frac{1}{r} = 1.
$$
If $P_i$ is superuniform under $H_i$, then $f_{\operatorname{BY}}(P_i)$ is an e-value for $H_i$ with $\EE[f_{\operatorname{BY}} (P_i)] \le 1$. Because $f_{\operatorname{BY}}$ is nonincreasing, the e-BH thresholds on the calibrated evidence coincide with the BY thresholds on the original p-values. 

For p-values $P = (P_1, P_2, \dots, P_K)$ and a uniform $U$ independent of $P$, the BY procedure is $D_{\operatorname{BY}} (P) = D_{\operatorname e}(f_{\operatorname{BY}} (P))$, its shared-uniform version is $\Uby (P, U) = D_U (f_{\operatorname{BY}} (P), U)$, and the optimal certified rule applied to the same calibrated evidence is
$$
\DstarBY (P, U) = D^\star (f_{\operatorname{BY}} (P), U) = \left\{i \in [K] : U \le \frac{\alpha}{K} \sum_{j = 1}^K \min\{f_{\operatorname{BY}} (P_i), f_{\operatorname{BY}} (P_j)\} \right\}.
$$

\begin{proposition} \label{prop-dstar-by}
Let $P_1, P_2, \dots, P_K$ be arbitrarily dependent, marginally superuniform p-values, and let $U \sim \operatorname{Uniform} (0, 1)$ be independent of $P$. Then $\DstarBY$ controls the FDR at $\alpha K_0 / K$, and
$$
D_{\operatorname{BY}} (P) \subseteq \Uby(P, U) \subseteq \DstarBY(P, U)
$$
for every realization of $U$.
\end{proposition}

\begin{proof}
Let $E_i = f_{\operatorname{BY}} (P_i)$, so that $\DstarBY (P, U) = D^\star (E, U)$, $D_{\operatorname{BY}} (P) = D_{\operatorname e} (E)$, and $\Uby (P, U) = D_U(E, U)$ by the identifications above. Theorem~\ref{thm-certification} then gives $\EE_U [Z_i(D^\star (E, U))] \le E_i$ for every $i$, and Proposition~\ref{prop-fdr-null} gives FDR control at $\alpha K_0 / K$. Since $D_{\operatorname e} (E) \subseteq D_U (E, U) \subseteq D^\star (E, U)$ by Proposition~\ref{prop-containment}, the same containments hold for $D_{\operatorname{BY}} (P)$, $\Uby (P, U)$, and $\DstarBY (P, U)$.
\end{proof}

\begin{figure}[!ht]
\centering
\begin{subfigure}{0.49 \textwidth}
\includegraphics[width = \textwidth]{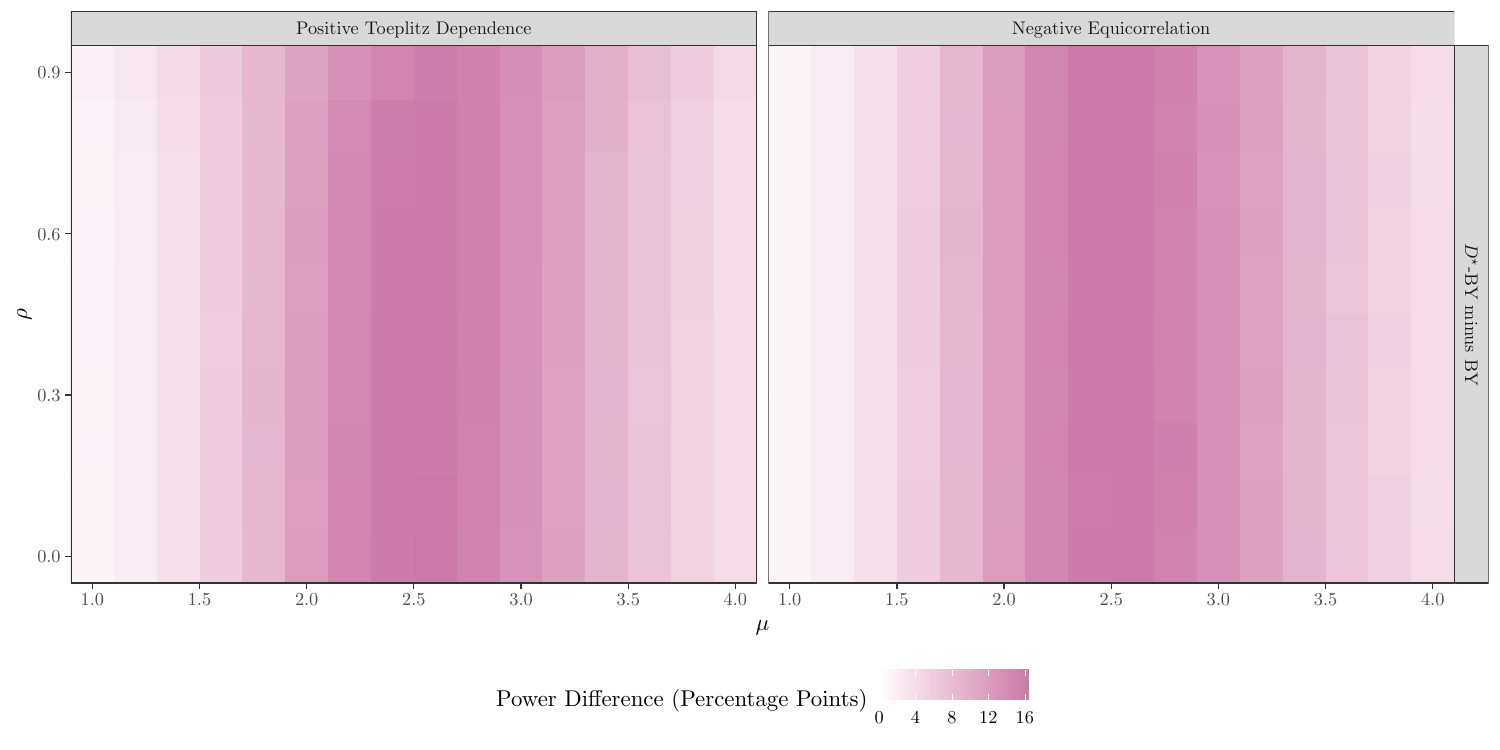}
\caption{$\DstarBY$ versus \BY}
\end{subfigure}\hfill
\begin{subfigure}{0.49 \textwidth}
\includegraphics[width = \textwidth]{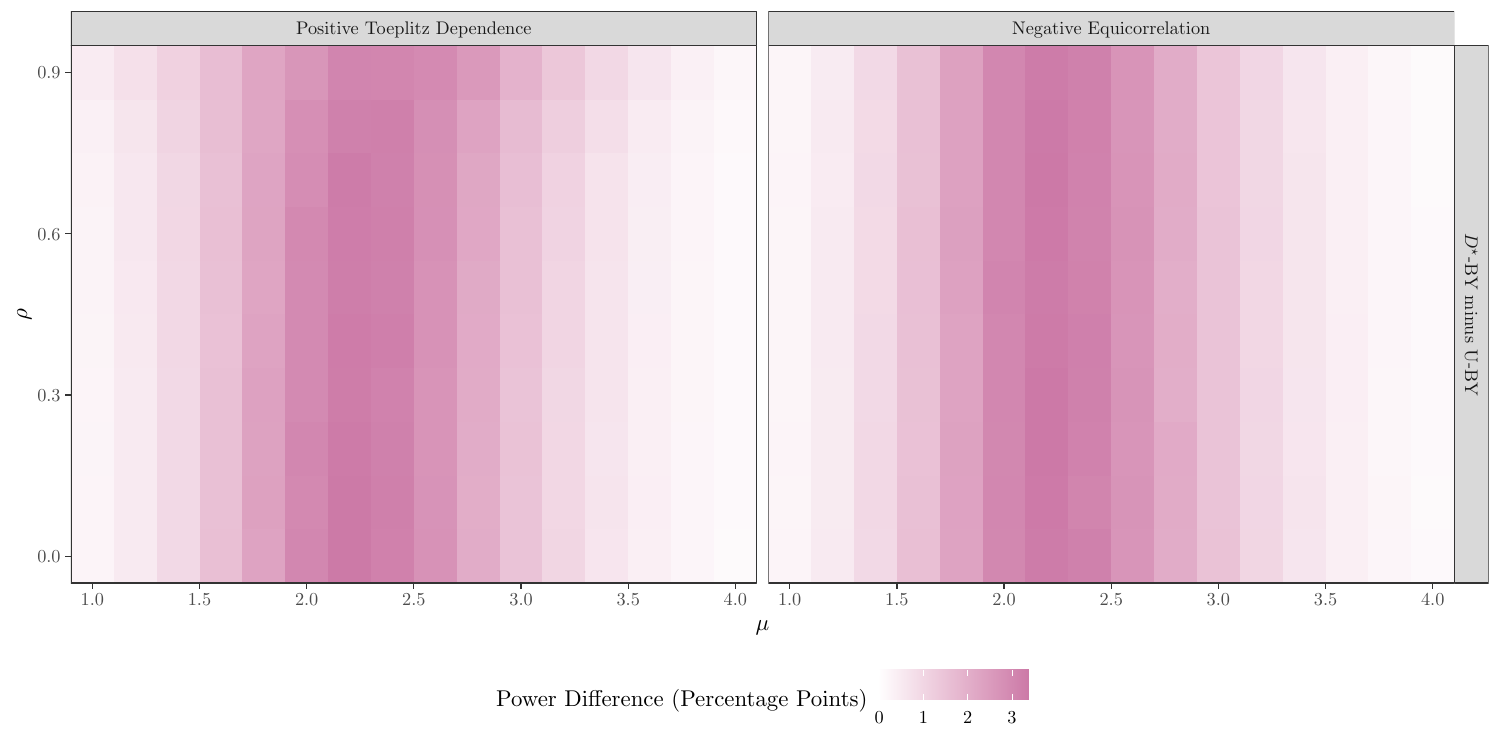}
\caption{$\DstarBY$ versus $\Uby$}
\end{subfigure}
\caption{Paired power differences, in percentage points, over the $(\mu, \rho)$ grid. Color intensities are not comparable across panels.}
\label{fig-supp-by-heatmaps}
\end{figure}

The simulation setup of Section~\ref{sec-simulation-setup} is used unchanged, with the BY-calibrated procedures added to the comparison and with the evidence supplied to the three procedures taken as $f_{\operatorname{BY}} (P_i)$ for $P_i = 1 - \Phi (Z_i)$, where $\Phi(\cdot)$ denotes the Gaussian CDF.

\begin{figure}[!ht]
\centering
\includegraphics[width = \textwidth]{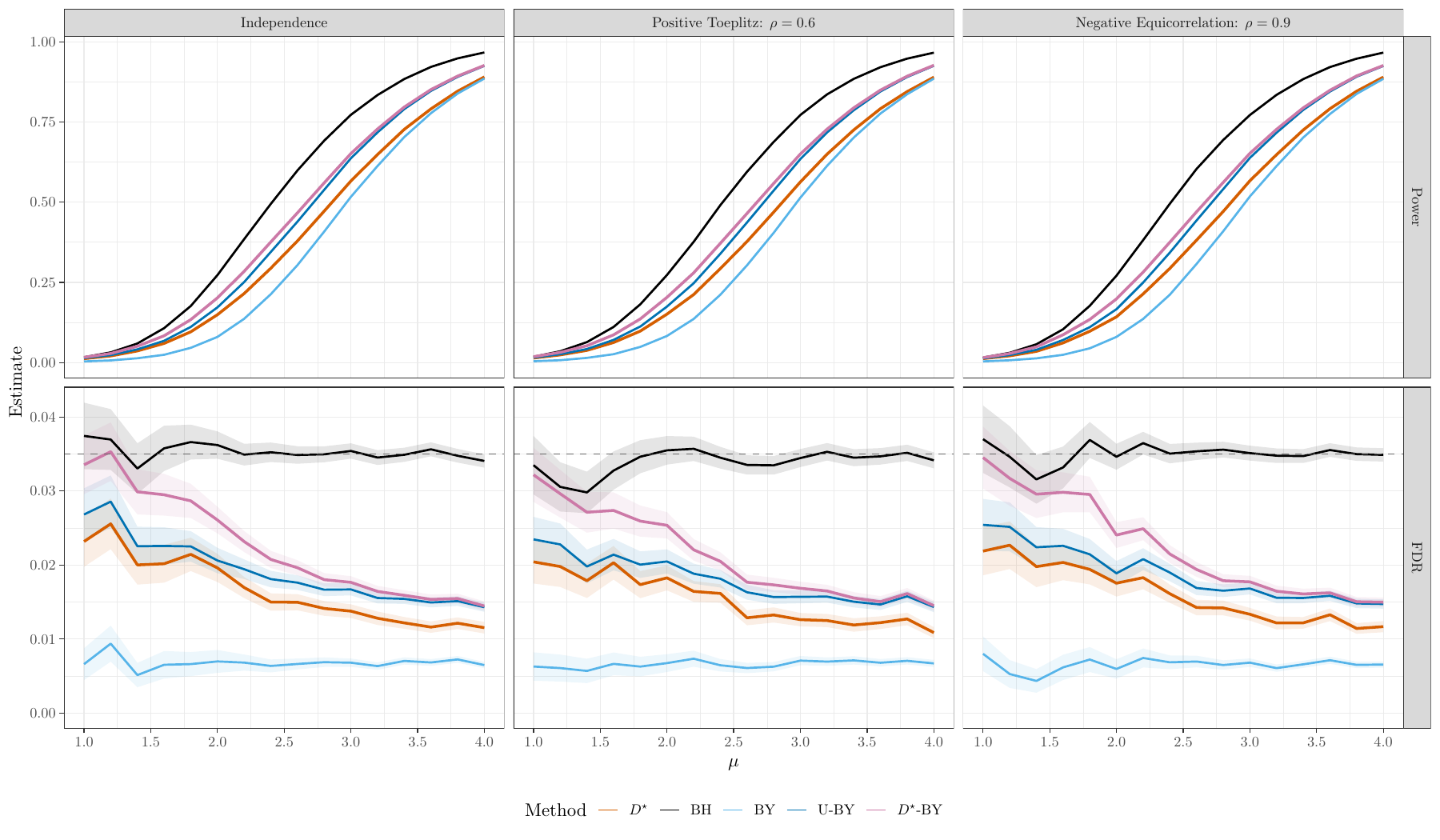}
\caption{Power (top row) and FDR (bottom row) against $\mu$ at independence ($\rho = 0$), positive Toeplitz dependence ($\rho = 0.6$), and negative equicorrelation ($\rho = 0.9$), for the BY-calibrated procedures and the \BH{} reference. Ribbons are pointwise $95\%$ Monte Carlo intervals. The dashed line in the lower panels marks $\alpha K_0 / K = 0.035$.}
\label{fig-supp-pvalue-benchmarks}
\end{figure}

Figure~\ref{fig-supp-by-heatmaps} reports the paired power differences for $\DstarBY$ against \BY{} and \Uby{}. The estimates reproduce the ordering $D_{\operatorname{BY}} \subseteq \Uby \subseteq \DstarBY$ of Proposition~\ref{prop-dstar-by}: $\DstarBY$ has higher power than both, and the gain is larger against BY than against \Uby{} at every setting. The gain concentrates at intermediate signal strengths in both comparisons and decays toward the endpoints of $\mu$.

Figure~\ref{fig-supp-pvalue-benchmarks} reports the power and FDR of the four  procedures, with Benjamini-Hochberg (\BH{}) procedure as a reference. The ordering $D_{\operatorname{BY}} \subseteq \Uby \subseteq \DstarBY$ of Proposition~\ref{prop-dstar-by} appears as a power ordering in every panel, and the FDRs of the three procedures remain at or below $\alpha K_0 / K$ throughout. \BH{} has higher power but no FDR guarantee in the negative-equicorrelation family.

\end{document}